\documentclass[10pt,a4paper,draft]{article}
\usepackage{amsmath,amssymb,amsthm,amscd}
\let\excellent\relax

\newcommand\al{\alpha}
\newcommand\be{\beta}
\newcommand\eps{\varepsilon}
\newcommand\la{\lambda}
\newcommand\ph{\varphi}

\newcommand\bU{\mathbf U}
\newcommand\BA{\mathbb A}

\newcommand\BP{\mathbb P}

\newcommand\C{{\cal C}}
\newcommand\D{{\cal D}}
\newcommand\CF{{\cal F}}

\newcommand\CO{{\cal O}}
\newcommand\CR{{\cal R}}
\newcommand\CX{{\cal X}}
\newcommand\CY{{\cal Y}}
\newcommand\CZ{{\cal Z}}

\newcommand\gm{\mathfrak m}

\newcommand\gp{{\mathfrak p}}

\DeclareMathOperator\Nil{Nil} \DeclareMathOperator\Spec{Spec}
\DeclareMathOperator\gr{gr} 
\DeclareMathOperator\Ker{Ker} \DeclareMathOperator\Ima{Im}
 
 \DeclareMathOperator\di{div}
\DeclareMathOperator\Supp{Supp}

\newcommand\isom{\simeq}

\newcommand\chr{\operatorname{char}}
\newcommand\Gal{\operatorname{Gal}}
\newcommand\he{\operatorname{height}}

 \DeclareMathOperator\Sw{{Sw}}

\newcommand\term[1]{{\it#1}}

\theoremstyle{plain} \swapnumbers
\newtheorem{pr}[subsection]{Proposition}
\newtheorem{thm}[subsection]{Theorem}
\newtheorem{lm}[subsection]{Lemma}
\newtheorem{cor}{Corollary}[subsection]

\newtheorem{quest}[subsection]{Question}

\theoremstyle{definition}

\newtheorem{rk}[cor]{Remark}
\newtheorem*{pf}{Proof}

\newcommand\paper[1]{\emph{#1}, }
\newcommand\jour[1]{#1 }
\newcommand\vol[1]{{\bf#1}}
\newcommand\yr[1]{ (#1)}

\begin{document}

\author{Igor Zhukov}
\title{Ramification of surfaces: sufficient jet order for wild jumps}

\maketitle

There exist different approaches to ramification theory of $n$-dimensional
schemes with $n\ge2$. One of these approaches is based on the following idea:
to reduce the situation of an $n$-dimensional scheme $\CX$ to a number of
$1$-dimensional settings, just by restricting to curves $C$
properly crossing the ramification subscheme $R\subset\CX$.
This way seems to be very natural; however, it has not got much attention,
except for \cite{Del:l} and \cite{Br:sur}. With the present paper
we intend to initiate the systematic development of this idea.

Let us be more precise. Let $L/K$ be a finite Galois extension of
the function field of a connected normal $n$-dimensional scheme
$\CX$, and let $f:\CY\to\CX$ be the normalization of $\CX$ in $L$.
Denote by $R\subset\CX$ the reduced branch locus of this morphism.
Let $C$ be a curve, i.~e., a closed integral 1-dimensional
subscheme of $\CX$. Since $f$ is finite, there are several curves
on $\CY$ lying over $C$; let $D$ be any of them. Then the natural
morphism $D\to C$ has familiar ramification invariants that can be
non-vanishing only at the points where $C$ meets $R$. The
principle is to collect these invariants for all regular curves on
$C$ that are distinct from the components of $R$ and to consider
the data obtained this way as a system of invariants of $f$.

One of the distant goals is to develop a sufficiently rich system
of ramification invariants that enables one to express numerous
known ramification data in terms of this system. On the other
hand, we hope to obtain a system of invariants with nice
functorial properties. For the latter reason, we do not require
$C$ to be transversal to $R$, and this is a new instant in this
paper. It seems that the consideration of curves tangent to the
branch divisor gives deep information; this can be observed
from consideration of Artin-Schreier extensions and their composites.
(Note that for 1-dimensional schemes, Swan characters at
ramified points form a system of invariants that enjoy both
properties: it is sufficiently rich and functorial.)

In the present paper we restrict the setting as follows:

$\bullet$ $n=2$;

$\bullet$ $\CX$ is equicharacteristic, i.~e., it can be equipped with
a morphism to $\Spec k$, where $k$ is an algebraically closed field
of prime characteristic $p$;

$\bullet$ If $P$ is any closed point of $\CX$, the local ring $\CO_{\CX,P}$  is
a 2-dimensional excellent local ring with the residue
field $k$;

$\bullet$ $\CX$ is regular;

$\bullet$ $L/K$ is solvable.

In section \ref{cal} we state a series of questions related to the
behavior of ramification jumps as one varies the curve $C$.
Some of these questions are new and some generalize the questions
considered in the above-mentioned papers.

In the following sections we answer the first of these questions
affirmatively; this is the main result of the paper.
Namely, we prove that the wild ramification jumps of the extension of the function
field of $C$ determined by $L/K$ at a point $P\in C\cap R$ depend only
on the jet of $C$ at $P$ of certain order, and we can bound this order
in a certain uniform way (see Corollary \ref{C1}).
Here $C\subset\CX$ is a curve which is regular at $P$.
Roughly speaking, if two curves $C$ and $C'$ have sufficiently high order
of tangence at $P$, then the wild jumps for $C$ and $C'$
are the same. By the wild jumps we mean the usual ramification
jumps (in lower numbering) divided by the index of tame ramification.

The proof of the main result is based on the work with arcs (parameterized algebroid
curves) on $\CX$, including the singular ones. One of the central ingredients
is a comparison of invariants of singularity of an arc on $\CX$
and those of an arc on $\CY$ above it in the case of cyclic $L/K$ of prime
degree; see Propositions \ref{lM} and \ref{pM}.

It would be important to prove a stronger fact: not only the
jumps but the whole ramification filtration is the same for $C$
and $C'$.
This statement (for curves transversal to $R$) is a step in Deligne's
program \cite{Del:l} describing how to compute Euler-Poincar\'e
characteristics of constructible \'etale sheaves on surfaces.

At present, we are able to prove the above fact for some class
of Galois groups that includes, in particular, all abelian $p$-groups, see Theorem
\ref{T2}.

As for the other questions from section \ref{cal}, in the case of
Artin-Schreier extensions
we can answer most of them affirmatively. This is the subject
of another paper \cite{ram-AS}.

\smallskip

I would like to thank P.~Deligne, L.~Illusie, G.~Laumon,
A.~N.~Parshin for the discussions that gave me inspiration for
this work and for the whole program. I am sincerely grateful to
I.~B.~Fesenko, B.~K\"ock, T.~Saito for numerous discussions and
also for the fact that they found errors in earlier versions of
the paper and gave me many suggestions on the improvement of
exposition. I would like to thank M.~V.~Bondarko, N.~L.~Gordeev
and I.~A.~Panin who answered my questions in commutative algebra
and algebraic geometry that arose in the course of this work. I am
deeply thankful to V.~P.~Snaith for the invitation to Southampton,
where a big part of the work was done, as well as for showing me a
lot of beautiful mathematics and for his all-embracing help during
my visits to Southampton.

I acknowledge the hospitality of Southampton University and
of Institut des Hautes \'Etudes Scientifiques (Bures-sur-Yvette)
and the financial support from
Royal Society and from RFBR (projects 00-01-00140 and
01-01-00997.)

\section{Terminology, notation, preliminary facts\label{nota}}

\subsection*{General notation}
$k$ is always an algebraically closed field of characteristic
$p>0$.

For any commutative ring $A$ denote
$$
\Spec_1A=\{\gp\in\Spec A|\he\gp=1\}.
$$

We denote by $v_X$ the valuation in the discrete valuation ring
$k[[X]]$.

For a prime ideal $\gp\in\Spec_1A$,
denote by $F_\gp$ the prime divisor $\Spec(A/\gp)$. If $\gp$ is a
principal ideal $(t)$, we write $F_t$ instead of $F_{(t)}$.

Let $A$ be an equal characteristic regular 2-dimensional local ring
with the maximal ideal $\gm=(T,U)$ and the residue field $k$,
$N$ a positive integer.
We have a canonical isomorphism between the completion of $A$
and $k[[T,U]]$.
A map $\lambda:k=A/\gm\to A$ is said to
be a section of level $N$ if the diagram
$$
\CD
k & @>\lambda>> & A
\\
@VVV && @VVV
\\
k[[T,U]]/(T,U)^{N}  & @= & A/\gm^{N}
\endCD
$$
commutes.

For any local ring $A$, we denote the completion of $A$ by
$\widehat{A}$.

The field of functions on an integral scheme $S$ is
denoted by $k(S)$.

\subsection*{Invariants of wild ramification}

Let $K$ be a complete discrete valuation field with perfect
residue field of characteristic $p$, and let $L/K$ be a finite Galois extension.
Let $v$ be the valuation on $L$ and $\CO_L=\{a\in L|v_L(a)\ge0\}$.

Recall that for any integral $i\ge-1$
the $i$th ramification subgroup in the group $\Gal(L/K)$ is defined as
\[
G_i=G_i(L/K)=\{g\in\Gal(L/K)|\forall a\in\CO_L:v(g(a)-a)\ge i+1\}.
\]
We have $|G_0|=e_{L/K}$, the ramification index, whereas
$|G_1|$ is the wild ramification index, i.~e., the maximal
power of $p$ dividing $e_{L/K}$. One can also define the tame
ramification index
\[
e^t(L/K)=(G_0:G_1).
\]

Now we introduce the ramification jumps, numbered from the bigger
to the smaller ones. Namely, for any $i\ge1$, denote
\[
h^{(i)}=h^{(i)}(L/K)=
\min\{j\ge1|p^i\text{ does not divide }|G_{j}|\}-1.
\]
This definition takes into account the multiplicities of jumps.
For example, if $G_{a}=G_{1}$ is of order $p^2$ for some
$a\ge1$ and $G_{j}$ is trivial for any $j>a$, then
$h^{(1)}=h^{(2)}=a$, and $h^{(i)}=0$ for $i\ge3$.
We define also modified (``wild'') ramification jumps
$w^{(i)}(L/K)=h^{(i)}(L/K)/e^t(L/K)$, $i=1,2,\dots$.
Note that $w^{(i)}(L/K)$ are non-negative and integral
in view of \cite[Ch. IV, Cor. 1 to Prop. 9]{CL}.

Next, we define the modified  Hasse-Herbrand
function $W_{L/K}:[1,\infty)\to[1,\infty)$ as
\[
W_{L/K}(u)=\int_1^u \frac {dt} {(G_1:G_t)}=\frac{\varphi_{L/K}(u)}{e^t(L/K)},
\]
where $G_u=G_i$, $i$ is the minimal integer such that $i\ge u$. Here
$\varphi_{L/K}$ is the usual Hasse-Herbrand function as defined
in \cite[Ch. IV, \S3]{CL}.

\cite[Ch. IV, Prop. 15]{CL} implies that for a Galois subextension
$M/K$ in $L/K$ we have $W_{L/K}=W_{M/K}\circ W_{L/M}$.

The following obvious property justifies the introduction
of modified jumps and modified Hasse-Herbrand function.

\begin{lm}
Let $K'/K$ be a tamely ramified extension.
Then $w^{(i)}(K'L/K')=w^{(i)}(L/K)$ for any $i$, and
$W_{K'L/K'}=W_{L/K}$.
\end{lm}

\subsection*{Ramification of surfaces}

A {\it surface} over $k$ is a connected normal \excellent 2-dimensional scheme
$\CX$ with a morphism $\CX\to\Spec k$ which induces an isomorphism
of residue fields at all points of codimension 2 in $\CX$.

Let $C,C'$ be distinct prime divisors on a surface $\CX$, and let $P$ be a closed
point of $\CX$. We define the intersection multiplicity of $C$ and $C'$ at $P$ as
\[
(C.C')_P=\begin{cases}\dim_k \CO_{\CX,P}/(\gp+\gp'), & P\in C\cap C',
\\
0, & P\notin C\cap C',
\end{cases}
\]
where $\gp$ and $\gp'$ are the prime ideals of $\CO_{\CX,P}$
corresponding to $C$ and $C'$.
By linearity this definition can be extended to any two
divisors $C,C'$ with no common components.

Let $\CX$ be a surface with the function field $K$.
Let $L/K$ be a finite separable extension, $\CY$ the normalization
of $\CX$ in $L$. We denote by $R_{L/K,\CX}\subset\CX$ the reduced
branch locus of the corresponding finite morphism $\ph:\CY\to\CX$,
i.~e., $R_{L/K,\CX}=\ph(\Supp\Omega^1_{\CY/\CX})$.
If $\CX$ is regular, then $R$ is of pure dimension 1. This follows
from the theorem of the purity of the branch locus in
\cite[\S41]{N}.

The extension $L/K$ is said to be \term{tame} (with respect
to $\CX$) if $L/K$ is tamely ramified with respect to any
extension of any discrete valuation of $K$ associated with
a prime divisor of $\CX$. Similarly one defines \term{unramified}
extensions.

A monoidal transformation means blowing up of a closed point.
The following fact is not difficult to prove. (Note that $R_{L/K,\CX_n}$ is a
part of the total transform of $R_{L/K,\CX}$.) This proposition
is believed to be standard, though we have not found a suitable reference.

\begin{pr} \label{NC}
Let $\CX$ be a regular surface over $k$, $L/K$ a finite extension
of its fraction field. Then there exists a sequence of monoidal
transformations $\CX_n\to\CX$ such that $R_{L/K,\CX_n}$ is
a simply normal crossing divisor.
\end{pr}

From now on, let $L/K$ be also normal.

Let $C$ be a prime divisor on $\CX$ such that $C\not\subset R_{L/K,\CX}$.
We can introduce $D_{C_1}$, the decomposition subgroup in $\Gal(L/K)$
at $C_1$, where $C_1$ is a component of $\ph^{-1}(C)$.
Since $L/K$ was not ramified at the generic point of $C$,
the group $D_{C_1}$ can be identified with $\Gal(k(C_1)/k(C))$.

Next, let $P$ be a regular point on $C$. Let $v$ be the corresponding
valuation on $k(C)$. Fix any extension $v_1$ of $v$ onto
the field $k(C_1)$. Let $K_{C,v}$ be the completion
of $k(C)$ with respect to $v$, and let $L_{C_1,v_1}$ be the completion
of $k(C_1)$ with respect to $v_1$. Then the decomposition group
$D_{C_1,v_1}$ of $v_1$ can be identified with $\Gal(L_{C_1,v_1}/K_{C,v})$.

Now we introduce
\[
w_{C,P}^{(i)}=w_{C,P}^{(i)}(L/K)=w^{(i)}(L_{C_1,v_1}/K_{C,v})
\]
and
\[
W_{C,P,L/K}=W_{L_{C_1,v_1}/K_{C,v}}.
\]
These objects do not depend on the choice of $C_1$ over $C$ and of $v_1$
since $\Gal(L/K)$ acts transitively on the set of such $C_1$,
and $\Gal(k(C_1)/k(C))$ acts transitively on the set of possible $v_1$.

Let $C$ and $C'$ be prime divisors on a surface $\CX$ and let $P$ be a closed point
in $C\cap C'$ such that both $C$ and $C'$ are regular at $P$.
We say that $L/K$ has \term{equal wild jumps with respect to $C$ and $C'$ at $P$},
if $w_{C,P}^{(i)}(L/K)=w_{C',P}^{(i)}(L/K)$ for any $i$.
We say that $L/K$ is \term{equally ramified with respect to $C$ and $C'$ at $P$
in the strong sense} if for any $C_1$ and $v_1$ one can choose
$C'_1$ and $v'_1$ such that $G_{j}(L_{C_1,v_1}/K_{C,v})=G_{j}(L_{C'_1,v'_1}/K_{C',v'})$
for any $j\ge-1$.

Let $M$ be a finite-dimensional representation of
$\Gal(L/K)$ (over any field $F$, $\chr F\ne p$).
Denote by $\Sw_{C,P}(M)$
the Swan conductor of the restriction of this representation on
the subgroup $D_{C_1,v_1}$.
In other words,
\begin{equation}
\label{sw-def}
\Sw_{C,P}(M)=\sum_{i=1}^\infty
 \frac{\dim_F(M/M^{G_{i}(L_{C_1,v_1}/K_{C,v})})}
 {(G_{0}(L_{C_1,v_1}/K_{C,v}):G_{i}(L_{C_1,v_1}/K_{C,v}))}.
\end{equation}
It is easy to see from this formula that $\Sw_{C,P}(M)$ does not
depend on the choice of $C_1$ over $C$ and of $v_1$.

\section{Questions\label{cal}}

Let $\CX$ and $L/K$ be as in the introduction.
Let $\CY$ be the normalization of $\CX$ in $L$.
Denote $R=R_{L/K,\CX}$.




For a closed point $P$ of $\CX$ denote
by $U_P$ the set of all prime divisors $C$ of $\CX$
such that $C$ is regular at $P$ and $C\not\subset R$.

\begin{quest}
\label{eusu}
\emph{(existence of a uniform sufficient jet order)}
Does there exists an effective divisor $R_0$ supported
at $R$ such that for any $P\in R$ and any $C,C'\in U_P$
the condition $(C.C')_P\ge (C.R_0)_P$ implies that
$L/K$ has equal wild jumps with respect to $C$ and $C'$ at $P$?
\end{quest}

If the answer is positive, we have $R_0=\sum_{i=1}^s m_iC_i$,
where $C_1\dots,C_s$ are all prime components of $R$.
Then we say that $m_i$ is a uniform sufficient jet order
at $C_i$, $i=1,\dots,s$.

To simplify the statements of the other questions, we make two technical
assumptions.

(i) $\CX$ is local, i.~e., $\CX=\Spec A$, where $A$ is a
2-dimensional local ring.
This enables us to write $(C.C')$ instead of $(C.C')_{(0)}$.
We shall also abbreviate $w^{(i)}_{F_\gp}$ as $w^{(i)}_{\gp}$.

(ii) $R$ is a simply normal crossing divisor in $\CX$.
(In particular, the components of $R$ are regular.)

Since $R$ is a simply normal crossing divisor, it consists
of at most 2 irreducible components. We shall not consider
the case $R=\emptyset$. If there is 1
component (resp. 2 components), we denote it by $C_1$
(resp. by $C_1$ and $C_2$); in the former case choose any regular
$C_2$ with $(C_1.C_2)=1$.

In both cases $(C_1.C_2)=1$; $C_1=F_{\gp_1}$, $C_2=F_{\gp_2}$;
we may introduce regular local
parameters $T$ and $U$ in $A$ such that $\gp_1=(T)$,
$\gp_2=(U)$.

Denote
\[
U_A=\{\gp\in\Spec_1 A|A/\gp\text{ is regular},\,F_\gp\not\subset R\}.
\]

Let $\gp\in U_A$. The set
$$
J_m(\gp)=\{\gp'\in U_A|(F_\gp.F_{\gp'})\ge m\}
$$
is said to be the jet of $\gp$ of order $m$. Also we introduce
$$
T_r=\{\gp\in U_A|(F_\gp.C_1)=r, (F_\gp.R)\le r+1\}
$$
and
$$
T_{r,m}=\{J_m(\gp)|\gp\in T_r\}.
$$

\rk We mentioned $R$ in the definition of $T_r$ in order
to exclude the curves tangent to $C_2$ from $T_1$ in the case
of two-component branch divisor.
\endrk

For the following question we have to endow the sets $T_{r,m}$
with the structure of affine varieties.

Fix a positive integer $n$ and a section $\lambda:A/\gm\to A$ of level $n$.

Let $\gp\in T_1$ (i.~e., this is the ideal of the germ of a curve
transversal to all components of $R$). Then $\gp=(f)$, where $f\equiv -U\mod(T,U^2)$.
Without loss of generality, we may assume
$$
f\equiv-U+\la(\al_1)T+\dots+\la(\al_n)T^n\mod\deg n+1,
$$
where $\al_1,\dots,\al_n\in k$ are determined uniquely by $\gp$.
(This follows from Weierstra\ss\ preparation theorem.)
Thus, if $R=C_1$, we can identify $T_{1,n}$ with the set of closed points of $\BA_k^n$
via:
$$
(\al_1,\dots,\al_n) \mapsto J_n((-U+\la(\al_1)T+\dots+\la(\al_n)T^n)).
$$
In the case $R=C_1+C_2$, the same map identifies $T_{1,n}$
with the set of the closed points of $(\BA^n_k)_{x_1\ne0}$î
Observe that $\al_1,\al_2,\dots$ are in fact the coefficients
in the expansion
$$
u=\al_1t+\al_2t^2+\dots,
$$
where $t$ and $u$ are the images of $T$ and $U$ in the completion
$\widehat{(A/\gp)}\simeq k[[t]]$. Therefore, $\al_1,\al_2,\dots$
are independent of $n$.

Similarly, if $\gp\in T_r$, $r\ge2$, we have $\gp=(f)$,
$$
f\equiv-T+\la(\be_r)U^r+\dots+\la(\be_n)U^n\mod\deg{n+1}.
\label{(2)}
$$
We have a bijection
\begin{align*}
(\BA_k^{n+1-r})_{x_1\ne0} & \to T_{r,n},
\\
(\be_r,\dots,\be_n) & \mapsto J_n((-T+\la(\be_r)U^r+\dots+\la(\be_n)U^n)).
\end{align*}
As in the previous case, $\be_1,\be_2,\dots$ are independent of $n$.

\begin{quest}
\emph{(semi-continuity of a jump)}
Let $m_1,m_2$ be uniform sufficient jet orders at $C_1,C_2$, fix any
positive integers $r$, $i$, and $m\ge rm_1+m_2$.
Does the set
$$\{J_m(\gp)|\gp\in T_r;\,w_\gp^{(i)}(L/K)\le s\}$$
form a closed subset in $T_{r,m}$ for any $s\ge0$?
\end{quest}

\begin{quest}
\emph{(generic value of a jump)}
\label{c:gen}
Fix any $r$
and $i$. Is
$$
w_r^{(i)}(L/K)=\sup\{w_\gp^{(i)}(L/K)|\gp\in T_r\}
$$
finite?
\end{quest}

\begin{quest}
\emph{(asymptotic of jumps)}
\label{c:as}
Is the sequence $(w_r^{(i)}(L/K)/r)_r$ (for any fixed $i$)
convergent and bounded by its limit?
\end{quest}

\begin{rk}
Consideration of Artin-Schreier extensions suggests affirmative answers
to all these questions, see \cite{ram-AS}.
\end{rk}

\subsection*{Representation version}

One can state similar questions where
$w_{C,P}^{(i)}(L/K)$ is replaced with $\Sw_{C,P}(M)$,
starting with

\begin{quest}
\label{eusu-sw}
\emph{(existence of a uniform sufficient jet order for Swan conductor)}
Let $M$ be a finite-dimensional representation of $\Gal(L/K)$ over any field
$F$, $\chr F\ne p$.
Does there exists an effective divisor $R_0$ supported
at $R$ such that for any $P\in R$ and any $C,C'\in U_P$
the condition $(C.C')_P\ge (C.R_0)_P$ implies
$\Sw_{C,P}(M)=\Sw_{C',P}(M)$?
\end{quest}

\begin{rk}
\label{suffstrong}
To answer this question affirmatively,
it is sufficient, in view of (\ref{sw-def}), to prove
the existence of uniform sufficient jet order for equal
ramification in the strong sense.
\end{rk}

\begin{rk} Versions of some of the questions,
in terms of Swan (or Artin) conductors, were
considered in \cite{Del:l} and \cite{Br:sur}.
Brylinski gives in \cite{Br:sur} the affirmative answer to
Question \ref{eusu-sw} for cyclic $p$-extensions under the
following assumptions:

$\bullet$ $\CX$ is a smooth algebraic surface;

$\bullet$ The branch divisor is a regular curve $D$;

$\bullet$ If a character of an extension corresponds
to the Witt vector $(x_0,\dots,x_{n-1})$, then all $x_i$
have neither poles nor zeroes  other than $D$ in the neighborhood of $P$
(it follows that $P$ is a non-exceptional point on the
branch divisor in Deligne's terminology);

$\bullet$ $(C.D)_P=1$ (i.~e., for curves transversal to the
branch divisor; Brylinski states that this condition
is not necessary).

In the case $(C.D)_P=1$ Brylinski computes the sufficient jet
order and the generic value of Artin conductor in terms
of Kato-Swan conductor.
\end{rk}

\begin{rk}
A representation version of Question \ref{c:as},
after extension to non-local $\CX$,
includes the existence of
\[
\Sw_\infty(M;D)=\lim_{r\to\infty}\frac{\Sw_r(M;D)}r,
\]
where $D$ is any prime component of $R$, and $\Sw_r(M;D)$ is
the generic value of $\Sw_{C,P}(M)$ over the $P\in D$ and
$C\in U_P$ with $(C.D)_P=r$.

We expect that $\Sw_\infty(M;D)$ is an important invariant
of ramification of $M$ at the place $D$, generalizing
Kato-Swan conductor, in a sense. However, at the moment we cannot
suggest a precise statement.
\end{rk}

\section{Arcs}
Here an arc is a parameterized algebroid curve drawn on a scheme.
Consideration of arcs is the main tool to answer Question \ref{eusu}.
We use \cite{Camp} as a convenient reference source for basic facts
about algebroid curves.

The terminology and notation are partially standard (see
\cite{Camp}) and partially customized for our needs.

\subsection*{Algebroid curves}

An (irreducible) \term{algebroid curve} over $k$ is a complete
Noetherian 1-dimensional local domain $\CR$ with a coefficient subfield
$k$. By abuse of language, $\Spec\CR$ is also said to be
an algebroid curve.

Let $\CR$ be an algebroid curve; the integral closure of $\CR$ is a
complete discrete valuation ring. Let $v$ be the corresponding
valuation.

Let $\gm$ be the maximal ideal of $\CR$.
An element $x\in\gm$ is said to be a \term{transversal parameter}
of $\CR$, if $x+\gm^2$ is not nilpotent in $\gr_\gm\CR$.
(According to \cite{Camp}, this is equivalent
to the following characterization property: $v(x)$ is minimal in $\gm$.)

The \term{multiplicity} of $\CR$ is the only $e=e(\CR)$
such that $\dim_k(\CR/\gm^n)=en+q$ for $n$ large enough.
It is known \cite[Th. 1.4.7]{Camp} that $e(\CR)=v(u)$
for any transversal parameter $u$.

Obviously, $e(\CR)=1$ iff $\CR$ is a regular ring.

We recall the definition of (strict) quadratic transform of an
algebroid curve from \cite{Camp}.

Let $\CR$ be an irreducible algebroid curve
over $k$, $\gm$ the maximal ideal of $\CR$. For $x\in\gm$,
$x\ne0$, denote
\[
\CR_x=\CR[x^{-1}z|z\in\gm].
\]
Now, the \term{quadratic transform} of $\CR$ is defined as
$\CR_1=\CR_x$ for any transversal parameter $x$ of $\CR$.
(It is known that $\CR_1$ is well defined.)

Obviously, $\CR_1\ne\CR$ unless $\gm$ is a principal ideal, i.~e.,
$\CR$ is regular.

Let $\overline{\CR}$ be the integral closure of $\CR$. Then
$\CR\subset\CR_1\subset\overline{\CR}$. In particular,
\begin{equation}
e(\CR_1)\le e(\CR).
\label{eqt}
\end{equation}
Since $\overline{\CR}$ is a Noetherian
$\CR$-module, the chain of quadratic transforms
\[
\CR\subset\CR_1\subset\CR_2=(\CR_1)_1\subset\dots
\]
stabilizes, i.~e., $\CR_i=\overline{\CR}$ for some $i$.

The minimal such $i$ is denoted by $M(\CR)$; we put
$M(\CR)=0$ for a regular $R$. Thus, if $M(\CR)\ge1$,
we have $M(\CR_1)=M(\CR)-1$.

\begin{lm}
\label{M}
If $\CR\subset\CR'\subset\overline{\CR}$, then
$M(\CR')\le M(\CR)$.
\end{lm}

This is obvious.

\subsection*{Arcs}

Let $\CX$ be any scheme over a field $k$. By an arc on $\CX$ we
mean a non-constant $k$-morphism $\C:\Spec k[[X]]\to\CX$ which
maps the closed point to a closed point. If $\CX=\Spec A$, where
$A$ is a local ring, an arc on $\CX$ can be identified with a local
homomorphism of $k$-algebras $f^\C:A\to k[[X]]$ such that
$\Ker f^\C$ is not the maximal ideal of $A$.

The \emph{center} of $\C$ is defined as $\C((X))$ and will be denoted by
$P_\C$. For a fixed closed point $P\in\CX$, the set of all arcs
such that $P_\C=P$ can be identified with
the set of arcs on $\Spec\CO_{\CX,P}$.
This enables us to
concentrate on the case of local schemes.

Let $P=P_\C$. Then $f^\C$ determines a $k$-algebra homomomorphism
$\widehat{f^\C}:\widehat{\CO_{\CX,P}}\to k[[X]]$ which corresponds
to an arc $\widehat{\C}$ on $\Spec\widehat{\CO_{\CX,P}}$.
It is easy to see that if $B$ is a $k$-subalgebra (with 1) of $k[[X]]$ and
$B$ is complete with respect to induced topology, then either
$B=k$, or $B$ contains $k[[t]]$ for some $t\in Xk[[X]]$, $t\ne0$.
Since $\C$ is non-constant, $\CR^\C=\Ima \widehat{f^\C}\ne k$.
Therefore, $\CR^\C$ contains $k[[t]]$ with $t$ as above.
It is easy to see that $\CR^\C$ is integral over $k[[t]]$,
whence $\dim\CR^\C=1$.
Thus, $\CR^\C$ is an algebroid curve over $k$.

The \term{support} of $\C$ is defined as the image of $\widehat{\C}$ and
will be denoted by $[\C]$. It is the closed subscheme of $\Spec\widehat{\CO_{\CX,P}}$
defined by $\Ker\widehat{f^\C}$ and isomorphic to $\Spec\CR^\C$.

Now let $\CX=\Spec A$; $A$ a two-dimensional Noetherian local ring admitting the
coefficient field $k$; $\gm$ the maximal ideal of $A$.

Let $x_1,\dots,x_n$ be a system of generators of $\gm$. Then we
can write $\widehat A=k[[X_1,\dots,X_n]]/I$ so that $x_i=X_i+I$,
$i=1,\dots,n$. To determine an arc on $\CX$, it is
necessary and sufficient to choose $f^\C_i=f^\C(x_i)\in X k[[X]]$
such that $b(f^\C_1,\dots,f^\C_n)=0$ for all $b\in I$
and $f^\C_i\ne0$ for some $i$.

We have $\CR^\C=k[[f_1^\C,\dots,f_n^\C]]$. We can decompose $f^\C$ as
\[
\CD A @> \alpha_A >> k[[X_1,\dots,X_n]]/I @> \beta_\C >> \CR^\C @>
\gamma_\C
>>
k[[t_\C]] @> \delta_\C >> k[[X]],
\endCD
\]
where $\alpha_A$ is the completion homomorphism,
and $\delta_\C\circ\gamma_\C\circ\beta_\C=\widehat{f^\C}$.
Further, $\beta_\C$ is surjective, $\gamma_\C$
embeds $\CR^\C$ into its integral closure in the field of fractions,
$\delta_\C$ is injective and makes $k[[X]]$ into a finite
$k[[t_\C]]$-algebra.

The \term{degree} of $\C$ is defined as
$d_\C=[k((X)):k((t_\C))]=v_X(t_\C)$.
A \term{primitive} arc is an arc of degree 1.
The \term{multiplicity} of $\C$ is defined as
\[
E_{\C}=e(\CR^\C)\cdot d_\C.
\]
Finally, we put $M_{\C}=M(\CR^\C)$. An arc is said to be
\term{regular} if $E_{\C}=1$. Obviously, $\C$ is regular iff it is
primitive and $M_{\C}=0$.

Two arcs $\C,\C'$ on a $k$-scheme $\CX$ are said to be
\term{weakly equivalent}, if $P_\C=P_{\C'}$, and $[\C]=[\C']$. Informally,
weakly equivalent arcs are merely different parameterizations of the
same algebroid curve on $\CX$. To make this formal, we
define a \term{parameterization} of an embedded algebroid curve
$\Spec\CR\to\CX$ as a composition
$\begin{CD}\Spec k[[X]]@>{\Spec\al}>>\Spec\CR\to\CX\end{CD}$,
where $\al$ is any injective local homomorphism of $k$-algebras $\CR\to k[[X]]$.
A parameterization is said to be primitive if it is a primitive arc.

Next, $\C$ and $\C'$ are said to be \term{equivalent}, if
$f^{\C'}=\lambda\circ f^\C$, where $\lambda$ is a $k$-algebra
automorphism of $k[[X]]$.
Obviously, any arc is weakly equivalent to a primitive one which
is defined uniquely up to equivalence.

Now we define the intersection multiplicity of two arcs
$\C$ and $\D$ on any surface $\CX$ such that
$P_\C$ is a  regular point on $\CX$.
This intersection
multiplicity is either a non-negative integer or $\infty$.

First, we put $(\C.\D)=0$  if $P_\C\ne P_\D$. In the remaining case we replace
$\CX$ with $\Spec\CO_{\CX,P_\C}$ and assume that $\CX=\Spec A$,
where $A$ is a two-dimensional regular local domain.

Assume that $\D$ is primitive.
Then we define
\[
(\C.\D)=v_X(\widehat{f^\C}(g_\D)),
\]
where $g_\D$ is a generator of $\Ker\beta_\D$.
The definition in
\cite[2.3.1]{Camp} is just a particular case of this one when
$\C$ is also primitive.

It is easy to see that $(\C.\D)$
is unchanged if $\C$ or $\D$ is replaced with an equivalent arc.

If $\D$ is not necessarily primitive, put
\[
(\C.\D)=d_\D\cdot(\C.\widetilde{\D}),
\]
where $\widetilde{\D}$ is defined by
$f^{\widetilde{\D}}=\delta_{\widetilde{\D}}\circ\gamma_\D\circ\beta_\D\circ\alpha_A$,
and $\delta_{\widetilde{\D}}$ maps $t_\D$ to $X$.
In view of the remark in the previous paragraph, $(\C.\D)$ is independent
of the choice of $t_\D$.
The relation to the intersection multiplicities of ``embedded plane
curves'' $(\CR^\C,\CR^\D)$ from \cite[2.3]{Camp} is as follows:
\[
(\C.\D)=d_\C {d_\D} (\CR^\C,\CR^\D).
\]

Next, we define $(\C.D)$ where $\C$ is an arc on $\CX$ as above,
and $D$ is a Weil divisor on $\CX$.
Assume first that $\CX=\Spec A$, where $A$ is a \emph{complete} local ring.
We require linearity on $D$ and assume, therefore, that $D$ is a
prime divisor: $D=\Spec A/\gp$. Note that $A/\gp$ is an algebroid curve.
Then by definition $(\C.D)=(\C.\D)$,
where $\D$ is any primitive parameterization of $D$.
(If $(\C.\D)=\infty$, $(\C.D)$ is assumed undefined.)

For general $\CX$, let $\CX'=\Spec\widehat{\CO_{\CX,P_\C}}$,
and let $f:\CX'\to\CX$ be the natural morphism.
Then by definition $(\C.D)=(\C_1.f^*D)$, where
$\C_1$ is the unique arc on $\CX'$ such that
$f\circ\C_1=\C$.

\begin{pr}
\label{2def}
Let $P$ be a regular point of codimension 2 on a surface $\CX$. Let $C$ and $C'$ be
any distinct irreducible curves on $\CX$ such that $P$ is a regular point
on both $C$ and $C'$, and let $D$ be a divisor on $\CX$ such that $C$
is not a component of $D$. Let $g_C$ and
$g_{C'}$ be local equations of $C$ and $C'$ at $P$. Denote
by $\C$ and $\C'$ any primitive parameterizations of $\widehat{\CO_{\CX,P}}/(g_C)$
and $\widehat{\CO_{\CX,P}}/(g_{C'})$ respectively.
Then $(C.C')_P=(\C.\C')$, and $(C.D)_P=(\C.D)$.
\end{pr}

\begin{pf}
None of $(C.C')_P$, $(\C.\C')$, $(C.D)_P$, $(\C.D)$ changes, if one
replaces $\CX$ with $\Spec\widehat{\CO_{\CX,P}}$, replacing
$C$, $C'$, and $D$ with their pullbacks. Therefore,
we may assume without loss of generality that $\CX=\Spec A$,
where $A$ is a complete 2-dimensional local ring with a coefficient
field $k$. We may also assume that $D$ is a prime divisor, and
the second equality to prove is reduced to the first one.
Since $C$ is regular, and $\C$ is primitive, $f^\C$ is surjective.
Therefore,
\begin{align*}
(\C.\C')&=v_X(f^\C(g_{\C'}))    \\
&=\dim_k k[[X]]/(f^\C(g_{\C'})) \\
&=\dim_k A/(\Ker f^\C + (g_{\C'})) \\
&=\dim_k A/(\Ker f^\C + \Ker f^{\C'}) \\
&=\dim_k A/(g_C,g_{C'}) \\
&=(C.C')_P.
\qed
\end{align*}
\end{pf}

\begin{lm}
\label{E}
Let $\C$ be an arc on $\Spec A$, where $A$ is a two-dimensional
regular local ring; let $T,U$ be local parameters of $A$.
Then

1. $E_\C=\min((\C.F_T),(\C.F_U))$.

2. There exist $\al,\be\in k$ such that
$(\C.F_{\al T+\be_U})>E_\C$; $[\al:\be]$ is a uniquely defined
point on $\BP^1_k$.
\end{lm}

\begin{pf}
1. Let $t=f^\C(T)$, $u=f^\C(U)$. We have $(\C.F_T)=v_X(t)$,
and $(\C.F_U)=v_X(u)$. It remains to note that either
$t$ or $u$ is a transversal parameter of $\CR^\C$, and $v_X=d_\C v$,
where $v$ is the valuation associated with the integral closure of
$\CR^\C$.

2. We have $v_X(\al t+\be u)>\min(v_X(t),v_X(u))$,
where $[\al:\be]$ is a uniquely defined
point on $\BP^1_k$.
\qed
\end{pf}

\subsection*{Monoidal transformations}

Consider an arc $\C:\Spec k[[X]]\to\CX$ on a surface $\CX$.
Assume that $O=P_\C$ is a regular point.
Let $\CX_1$ be the blowing up of $\CX$ at
the point $O$. By the second part of Lemma \ref{E},
$\C$ determines a unique
point $O_1$ in the exceptional divisor, and we denote
$A=\CO_{\CX,O}$,
$A_1=\CO_{\CX_1,O_1}$. We can write down a commutative diagram
\[
\CD
A @> \alpha_A >> & k[[X_1,X_2]] @> \beta_\C >> & \CR^\C @> \gamma_\C >>
& k[[t_\C]] @> \delta_\C >> & k[[X]]
\\
@VVV & @VVV & @VVV & \| & @. \|
\\
A_1 @> \alpha_{A_1} >> & k[[X'_1,X'_2]]/I_1 @> \beta_{\C_1} >>
& (\CR^\C)_1 @> \gamma_{\C_1} >>
& k[[t_\C]] @> \delta_{\C_1} >> & k[[X]],
\endCD
\]
where
$(\CR^\C)_1$ is the quadratic transform of $\CR^\C$,
$\beta_{\C_1}$ comes from the embedding
of strict quadratic transform of an embedded
algebroid curve $\beta_C$ into the formal quadratic transformation
of $\Spec k[[X_1,X_2]]$ in the direction defined by $O_1$,
see the discussion in \cite[pp. 35--38]{Camp}; $\gamma_{\C_1}$ is
induced by $\gamma_\C$; $\delta_{\C_1}=\delta_\C$. The composition
$\delta_{\C_1}\circ\gamma_{\C_1}\circ\beta_{\C_1}\circ\alpha_{A_1}$
determines an arc $\C_1$ on $\Spec A_1$ and on $\CX_1$.
This arc is said to be the \term{strict transform} of $\C$.
It follows from (\ref{eqt}) that $E_{\C_1}\le E_{\C}$.

If $P_\C=P_\D$, then it is immediate from \cite[2.3.2]{Camp} that
\begin{equation}
(\C.\D)=(\D.\C)=E_{\C}\cdot E_{\D}+(\C_1.\D_1),
\label{iqt}
\end{equation}
where $\C_1$ and $\D_1$ are the strict transforms of $\C$ and $\D$
after blowing up $\CX$.

The point $O_1\in\CX_1$ is said to be \term{the first infinitely near
point} for $\C$. The \term{$(i+1)$th infinitely near point} for $\C$ is
defined as the $i$th infinitely near point for $\C_1$.

The lemma below follows easily from the definitions.

\begin{lm}
\label{EE} Let $\C$ be an arc such that $P_\C$ is a regular point.
Then $E_\C=(\C_1.E)$, where $\C_1$ is the strict transform of
$\C$, and $E$ is the exceptional divisor.
\end{lm}

\subsection*{Hamburger-Noether expansion}

Let $\C$ be an arc on $\Spec A$, where $A$ is a two-dimensional
regular local ring; let $T,U$ be local parameters of $A$.
Let $t=f^\C(T)$ and $u=f^\C(U)$.
Assume $v_X(t)\ge v_X(u)$.
Then one can write, in a unique way,
\begin{equation}
\label{HN}
\begin{aligned}
t = a_{01} u +& a_{02} u^2 + \dots + a_{0h_0}u^{h_0} + u^{h_0}z_1, \\
u = \qquad & a_{12} z_1^2 + \dots + a_{1h_1}z_1^{h_1} + z_1^{h_1}z_2, \\
& \dots \\
z_{r-2} = \qquad & a_{r-1,2} z_{r-1}^2  + \dots
        + a_{r-1,h_{r-1}}z_{r-1}^{h_{r-1}} + z_{r-1}^{h_{r-1}}z_r, \\
z_{r-1} = \qquad & a_{r2}z_r^2 + a_{r3}z_r^3+\dots,
\end{aligned}
\end{equation}
where $r\ge0$, $h_0,\dots,h_{r-1}$ are positive integers,
$a_{ij}\in k$, $z_i\in k[[X]]$,
$E_{\C}=v(u)>v(z_1)>\dots>v(z_r)=d_{\C}$,
and $a_{rj}\ne0$ for some $j$ if $r>0$. This is the Hamburger-Noether
expansion of $\C$ in parameters $T,U$.

Conversely, any system (\ref{HN}) is the Hamburger-Noether
expansion of some arc $\C$ on $A$ in parameters $T,U$. The data
$(r;(h_i);(a_{ij}))$ correspond bijectively to the classes
of weakly equivalent arcs on $\Spec A$.

We have $r=0$ iff $e(\CR^\C)=1$.
In this case the Hamburger-Noether expansion is just an
\emph{equation} of $[\C]$ in parameters $T,U$.

If $r>0$, we have
$h_0+h_1+\dots+h_{r-1}=M_{\C}$.

Let $\C$, $\D$ be two arcs on $\Spec A$. Then the length of the
``common part'' of Hamburger-Noether expansions of $\C$ and $\D$
is equal to the number of common infinitely near points for $\C$
and $\D$.

A detailed discussion of these and related facts on Hamburger-Noether
expansions is given in \cite{Camp}.

\subsection*{Ramification invariants associated with arcs}

Let $\CX$ be a surface over $k$; $L/K$ a finite Galois extension
of its function field; $\CY$ the normalization of $\CX$ in $L$;
$R=R_{L/K,\CX}$.

Let  $\bU_\CX=\bU_{\CX,L/K}$ consist of all arcs $\C$ on $\CX$ such that
$(\C.R)$ is defined, i.~e., $P_\C$ is a regular point,
and $[\C]$ is not a component of the pullback
of $R$ to $\Spec \widehat{\CO_{\CX,P_\C}}$. We shall define wild ramification jumps
$w^{(i)}_\C(L/K)$ for any $\C\in\bU_\CX$, $i=1,2,\dots$.

Let $\C\in\bU_\CX$. Replacing $\CX$ with
$\widehat{\CO_{\CX,P_\C}}$,
we may assume that $\CX=\Spec A$, where $A$ is a complete local
ring. Then $[\C]$ is a prime divisor on $\Spec A$. The prime divisors $D$ of $\CY$ lying
over $[\C]$ are associated with all extensions to $L$ of the
valuation on $K$ determined by $[\C]$. Therefore, $\Gal(L/K)$ acts
transitively on the set of all such $D$.

Fix any $D$, and let $H\subset\Gal(L/K)$ be the
stabilizer of $D$.
Then $H=\Gal(k(D)/k([\C]))$. Let $\CO$ and $\CO_1$
be the valuation rings in the complete discrete valuation
fields $k([\C])$ and $k(D)$ respectively.
Consider a co-Carthesian square
\[
\begin{CD}
k([\C]) & @>>> &  k((X))
\\
@VVV && @VVV
\\
k(D) & @>>> & k(D)\otimes_{k([\C])} k((X)),
\end{CD}
\]
where the upper arrow is induced by $\delta_\C$.
Then $k(D)\otimes_{k([\C])} k((X))$
 is an $H$-Galois algebra over $k((X))$.
It follows that $k(D)\otimes_{k([\C])} k((X))\isom\prod_{i=1}^dL_\C$ for some $d$
and some Galois extension
$L_\C/k((X))$ of degree $|H|/d$.

We introduce
\[
w_\C^{(i)}(L/K)=w^{(i)}(L_\C/k((X))), \quad i=1,2,\dots
\]
and
\[
W_{\C,L/K}=W_{L_\C/k((X))}.
\]

\begin{lm}
\label{ramarc}
Let $C$ be a prime divisor on $\CX$, $P$ a regular
closed point on $C$. Let $g\in\CO_{\CX,P}$ be a local equation of $C$
at $P$. Denote by $\C$ any primitive parameterization
of $\widehat{\CO_{\CX,P}}/(g)$. Then
$w_{C,P}^{(i)}(L/K)=w_\C^{(i)}(L/K)$.
\end{lm}

\begin{pf}
The extension $L_{C_1,v_1}/K_{C,v}$ from the definition
of $w_{C,P}^{(i)}(L/K)$ coincides with $k(D)/k([\C])$
from the definition of $w_\C^{(i)}(L/K)$. The latter extension
is isomorphic to $L_\C/k((X))$ since $\C$ is a primitive arc.
\qed
\end{pf}

\section{Main theorem}

Let $\CX$ be a surface over $k$; $L/K$ a finite solvable Galois extension
of its function field; $\CY$ the normalization of $\CX$ in $L$;
$R=R_{L/K,\CX}$.

We say that $L/K$ \term{has equal wild jumps at
$\C,\D\in\bU_{\CX,L/K}$}, if $w_\C^{(i)}(L/K)=w_\D^{(i)}(L/K)$ for all
$i=1,2,\dots$. This is equivalent to
$W_{\C,L/K}=W_{\D,L/K}$.

\begin{thm}
\label{main}
In the above setting, there exists a non-decreasing sequence
of positive integers $\Delta_\CX(L/K,i)$, $i=1,2,\dots$,
such that $\Delta_\CX(L/K,i)\ge i^2$, and if $\C,\D\in\bU_\CX$,
$P_\C=P_\D$ is a regular point on $\CX$,
$\C,\D$ are primitive, and
\[
(\C.\D)\ge ((\C.R)+\max(M_{\C},M_{\D}))\Delta_\CX(L/K,\max(E_{\C},E_{\D})),
\]
then $L/K$ has equal wild jumps at $\C$ and $\D$.
\end{thm}

Theorem \ref{main} will be proved in the following 4 sections.

\begin{cor}
\label{C1}
The answer to Question \ref{eusu} is affirmative.
\end{cor}

\begin{pf}
It is sufficient to put
\[
R_0=\Delta_\CX(L/K,1)R
\]
and to apply Prop. \ref{2def} and Lemma \ref{ramarc}.
\qed
\end{pf}

\begin{rk}
If $C$ is a regular curve on a surface $\CX$, and $\CY$
is the normalization of $\CX$ in a finite extension of its
function field, then the irreducible curves on $\CY$ over $C$
(``liftings'' of $C$) are, in general, not regular.
This is why we have to introduce $\Delta_\CX(L/K,i)$ with $i>1$ and to prove
Theorem \ref{main} which is more general than Corollary \ref{C1}.
\end{rk}
%

\section{Lifting of arcs}

\setcounter{cor}0

Let $f:\CY\to\CX$ be a finite surjective morphism of $k$-schemes,
$\C$ an arc on $\CX$. An arc $\C'$ on $\CY$ is said to be
a \term{lifting} of $\C$ onto $\CY$ if $\C':\Spec k[[X]]\to\CY$ can be
factored as $\C'=\C_f\circ g$, where $\C_f$ is determined by
a Carthesian square
\begin{equation*}
\begin{CD}
\Spec k[[X]]\times_{\CX}\CY & @>\C_f>> & \CY
\\
@VVV && @VfVV
\\
\Spec k[[X]] & @> \C >> & \CX
\end{CD}
\end{equation*}
and $g$ is the normalization of an irreducible component of $\Spec
k[[X]]\times_{\CX}\CY$ with reduced scheme structure.

\begin{rk}
\label{L1}
Let $\C$ be a primitive arc.  Analyzing a more
detailed diagram
\begin{equation*}
\begin{CD}
\Spec k[[X]]\times_{\CX}\CY & @>>> & [\C]\times_{\CX}\CY &
@>>> & \Spec\widehat{\CO_{\CX,P_\C}}\times_{\CX}\CY & @>>> & \CY
\\
@VVV && @VVV && @VVV && @VfVV
\\
\Spec k[[X]] & @>\Spec\gamma_\C>> & [\C] &
@>\Spec\beta_\C >> & \Spec\widehat{\CO_{\CX,P_\C}}& @>>> & \CX
\end{CD}
\end{equation*}
we see that the following numbers coincide:

1) the number of non-equivalent liftings of $\C$ onto $\CY$;

2) the number of non-equivalent liftings of
$\widehat{\C}:\Spec k[[X]]\to\Spec\widehat{\CO_{\CX,P_\C}}$
onto $\Spec\widehat{\CO_{\CX,P_\C}}\times_{\CX}\CY$;

3) the number of irreducible components in $[\C]\times_{\CX}\CY$
(note that $[\C]\times_{\CX}\CY$ is the spectrum of a finite
$\CR^\C$-algebra);

4) the number of irreducible components in $\Spec k[[X]]\times_{\CX}\CY$.

(See also proof of Prop. \ref{tamelift} below.)
\end{rk}

\begin{rk}
\label{L2}
It is also clear from the above diagram that any arc $\D$
on $\CY$ is a lifting of the arc $f\circ\D$ on $\CX$.
Moreover, if $\D$ is primitive, and $\C$ is a primitive
and arc on $\CX$ weakly equivalent to $f\circ\D$,
then $\D$ is also a lifting of $\C$. Thus, any primitive
arc on $\CY$ is a lifting of a primitive arc on $\CX$.
\end{rk}

We examine liftings in the case of some special finite extensions
of two-dimensional local rings.

\begin{pr}
\label{tamelift}
Let $A=k[[T,U]]$, $A_1=k[[T_1,U]]$, and $A$ is embedded into $A_1$
by mapping $T$ to $\xi(T_1,U)\equiv T_1^l\mod(T_1^{l+1},T_1U)$, where
$l$ is a prime number. Assume also that there exists an $A$-automorphism
$\sigma$ of $A_1$ such that $\sigma^l=1$, and
$A_1^{\langle\sigma\rangle}=A$.

1. Let $\C$ be a primitive arc on $\Spec A$ with $[\C]\ne F_T$.
Then the number $n_\C$ of non-equivalent liftings of $\C$ onto $\Spec A_1$
 is either $1$ or $l$, and these liftings are primitive arcs.

2. If $\C'$ is a lifting of a primitive arc $\C$ as above, we have
\begin{align*}
(\C'.F_{T_1}) &= (\C.F_T)/n_\C, \\
(\C'.F_U) &= l(\C.F_U)/n_\C, \\
E_{\C'} &=\min((\C.F_T),l(\C.F_U))/n_\C.
\end{align*}
If $n_\C=l$, we have $M_{\C'}\le M_\C$.

3. Let $\C$ and $\D$ be primitive arcs on $\Spec A$ with $[\C],[\D]\ne F_T$.
Let $\C'$ and $\D'$ be their liftings. If $n_\C=1$ or $n_\D=1$, then
$(\C'.\D')=l(\C.\D)/(n_\C n_\D)$. If $n_\C=n_\D=l$, then
$\sum_{i=1}^{l}(\C_i.\D')=(\C.\D)$, where $\C_1,\dots,\C_l$ are
all pairwise non-equivalent liftings of $\C$.
\end{pr}

\begin{pf}
Let $\C$ be an arc on $\Spec A$ with $[\C]\ne F_T$.
It is clear from Remark \ref{L1} that $n_\C$ coincides
with the number of primes in $\Spec A_1$ over the generic
point of $[\C]$. By \cite[Ch. V, \S2, Th.~2]{Bourbaki},
$\langle\sigma\rangle$ acts transitively on the set
of such primes, whence $n_\C$ is either $1$ or $l$.
Denote by $g$ any generator of $\Ker f^\C$.

Let first $n_\C=1$. Since $[\C]\ne F_T$, there are no nilpotents
in $\CR^\C\otimes_A A_1$. The unique lifting $\C'$ of $\C$
corresponds to the bottom row of the diagram
\begin{equation}
\label{cdl}
\begin{CD}
A & @> \beta_\C >> & A/(g)=\CR^\C & @> \gamma_\C >> & k[[t_\C]] &=& k[[X]]
\\
@VVV && @VVV && @VVV
\\
A_1 &@>>>& A_1/(g)=\CR^\C\otimes_AA_1 &  @>\gamma_\C\otimes_AA_1 >> &
B @>>> B/\Nil B
\end{CD}
\end{equation}
extended to the right with the embedding of $B/\Nil B$
into its integral closure $B'\isom k[[X]]$.

Since $A'$ is a free $A$-module of rank $l$, it is clear
that $\CR^\C\otimes_AA_1$ is of finite $k$-codimension
in $B$ and in $B'$, whence $\C'$ is a primitive arc.
Next, it is obvious that $g_1=g(\xi(T_1,U),U)\in k[[T_1,U]]$
is an equation of $\C'$. By the definition
and Lemma \ref{E} we have
\begin{align*}
(\C'.F_{T_1}) &= v_X(g_1(0,X)) = v_X(g(0,X)) =   (\C.F_T)/n_\C, \\
(\C'.F_U) &=  v_X(g_1(X,0)) = l v_X(g(X,0)) =  l(\C.F_U)/n_\C, \\
E_{\C'} &= \min((\C.F_{T_1}),(\C.F_U))
=\min((\C.F_T),l(\C.F_U))/n_\C.
\end{align*}

Next, let $n_\C=l$. Then $g(\xi(T_1,U),U)=g_1\dots g_l$ in $A_1$,
and $(g_1),\dots,(g_l)$ are distinct prime ideals.
Consider $(A_1/(g_1\dots g_l))\otimes_{A/(g)}k((t_\C))$.
This is a $k((t_\C))$-algebra of dimension $l$ with $l$
minimal prime ideals; therefore, it is isomorphic to
the direct product of $l$ copies of $k((t_\C))$.
Applying $\otimes_{A/(g)}k((t_\C))$ to the right
square of the diagram \eqref{cdl}, we see that
non-equivalent liftings of $\C$ are, up to equivalence,
the primitive arcs $\C_1,\dots,\C_l$ with $\CR^{\C_i}=A_1/(g_i)$.

The group $\langle\sigma\rangle$ acts transitively on
$\{(g_1),\dots,(g_l)\}$, therefore, on $\{\C_1,\dots,\C_l\}$.
It follows that $(\C_i.F_{T_1})$
and $(\C_i.F_U)$ are independent of $i$. On the other hand,
\begin{align*}
\sum_{i=1}^l(\C_i.F_{T_1}) &= \sum_{i=1}^lv(g_i(0,X)) = v(g(0,X)) =   (\C.F_T), \\
\sum_{i=1}^l(\C_i.F_U) &=  \sum_{i=1}^lv(g_i(X,0)) = l v(g(X,0)) =  l(\C.F_U),
\end{align*}
and the desired formulae for $(\C'.F_{T_1})$,
 $(\C'.F_U)$, $E_{\C'}$ follow.

Since $A_1/(g_i)\isom\CR^{\C_i}$ is integral over
$A/(g)\isom\CR^{\C}$ with the same field of fractions, it follows
from Lemma \ref{M} that $M_{\C_i}\le M_\C$.

To prove the remaining assertion, we have to relate $f^\D$ and
$f^{\D'}$. Let $\psi$ be the map that makes the diagram
\[
\begin{CD}
A & @>f^\D>> & k[[X]] \\
@VVV && @V\psi VV \\
A_1 & @>f^{\D'}>> & k[[X]]
\end{CD}
\]
commutative.
If $n_\D=1$, $\psi$ is an embedding of $k[[X]]$ into a totally ramified extension
of degree $l$. In particular,
\begin{align*}
\psi(f^\D(T))&=\xi(f^{\D'}(T_1),f^{\D'}(U)), \\
\psi(f^\D(U))&=f^{\D'}(U).
\end{align*}
In this case $(\C'.\D')$ is computed exactly in the same way as
$(\C'.F_U)$.

If $n_\D=l$, $n_\C=1$, we change the roles of $\C$ and $\D$.

Finally, if $n_\C=n_\D=l$, then $\psi$ is the identity map, and
\begin{align*}
f^\D(T)&=\xi(f^{\D'}(T_1),f^{\D'}(U)), \\
f^\D(U)&=f^{\D'}(U).
\end{align*}
We obtain
\begin{align*}
\sum_{i=1}^{l}(\C_i.\D')&=\sum_{i=1}^{l}v_X(g_i(f^{\D'}(T_1),f^{\D'}(U)))
\\
&=v_X(g(\xi(f^{\D'}(T_1),f^{\D'}(U)),f^{\D'}(U)))
\\
&=v_X(g(f^\D(T),f^\D(U)))=(\C.\D).
\qed
\end{align*}
\end{pf}

\begin{cor}
In the setting of Proposition we have $E_{\C'}\le l E_{\C}$.
\end{cor}

\begin{pf}
Use Lemma \ref{E}.
\qed
\end{pf}

\subsection*{Growth of $M_{\C}$ in cyclic extensions of prime degree}

In this subsection $\CX$ is a regular $k$-surface with the function field
$K$; $L/K$ is a cyclic extension of prime degree $l$;
 $\CY\to\CX$ is the normalization of $\CX$ in $L$.
For any arc $\C$ on $\CX$ we denote by $n_\C$ the number
of non-equivalent liftings of $\C$ onto $\CY$.


Up to the end of the section we assume that $l\ne p$. Similar
statements with $l=p$ are deferred until the analysis of
Artin-Schreier coverings in the section \ref{AS}.

\begin{lm}
\label{nc}
Let $\C\in\bU_\CX$ be a primitive arc.
Assume that $L=K(x)$, $x^l=a$, where $a\in\CO_{\CX,P_\C}$.
Then $n_\C=1$ iff $l\nmid v_X(f^\C(a))$.
\end{lm}

\begin{pf}
In view of Remark \ref{L1}, we may assume $\CX=\Spec A$,
where $A$ is a local ring. We have $\CY=\Spec B$,
where $B$ is the integral closure of $A[x]$.
Let $U$ be the complement of $R_{L/K,\CX}$ in $\CX$.
The unramified morphism $\CY\times_\CX U\to U$
factors through $\Spec A[x]\times_\CX U$, whence
$\Spec A[x]\times_\CX U$ is regular, and
$\CY\times_\CX U=\Spec A[x]\times_\CX U$.

By Remark \ref{L1}, $n_\C$ coincides with
the number of irreducible components in $[\C]\times_\CX\CY$,
i.~e., the number of primes in $B\otimes_A\hat A$ over
the generic point $\gp$ of $[\C]$. Since $\C\in \bU_{\CX}$,
we obtain that
$\gp$ is in $\Spec\hat A\times_\CX U$, and $n_\C$ is the number of primes
in $A[x]\otimes_A\hat A$ over $\gp$, i.~e., the number of minimal
primes in
\[
(\hat A/\gp)[Y]/(Y^l-a)\isom
k[[X]][Y](Y^l-f^\C(a))\isom
k[[X]][Y](Y^l-X^m),
\]
where $m=f^\C(a)$. If $l|m$, then $Y^l-X^m$
is a product of $l$ irreducible factors and $n_\C=l$.
If $l\nmid m$, it is easy to check that $Y^l-X^m$
is irreducible, and $n_\C=1$. (Indeed, apply induction
on $\max(m,l)$, and substitute either $XZ$ for $Y$
or $YZ$ for $X$.)
\qed
\end{pf}

\begin{lm}
\label{l}
Assume that $R_{L/K,\CX}$ is a simply normal crossing divisor.
Let $\C\in\bU_\CX$ be a primitive arc with $n_\C=1$, not regular.
Let $\D\in\bU_\CX$ be a primitive arc such that $E_{\D}<l E_{\C}$,
and $(\C.\D)>(M_{\C}+1)E_{\C}E_{\D}$.
Then

1. $n_\D=1$;

2. $E_{\D}\ge E_{\C}$.
\end{lm}

\begin{pf}
Let $n=M_{\C}$. Then $\C$ and $\D$ have at least $n$ common
infinitely near points $O_1,\dots,O_n$. Let $\C_n$ and $\D_n$
be the $n$th strict transforms of $\C$ and $\D$ respectively.
Then $\C_n$ is regular and the only line of the Hamburger-Noether
expansion of $\C_n$ (in suitable local parameters $t$ and $u$ at $O_n$)
coincides with the last line of the
Hamburger-Noether expansion of $\C$ (in some local parameters $t_0,u_0$):
\[
z_{r-1}=a_{r2}z_r^2+a_{r3}z_r^3+\dots,
\]
where $z_{r-1}=f^{\C_n}(t)$, $z_r=f^{\C_n}(u)$.
Let $j=\min\{j|a_{rj}\ne0\}$. Then $j=v_X(z_{r-1})=E_{\C_{n-1}}\le E_{\C}$, and
\[
(\C_n.\D_n)\ge(\C.\D)-M_{\C}E_{\C}E_{\D}>E_{\C}E_{\D}\ge j
E_{\C_n}E_{\D_n}.
\]
Therefore, $\C_n$ and $\D_n$ have at least $j$ common infinitely
near points. It follows that the beginning of the first line
in the Hamburger-Noether expansion of $\D_n$ with respect to $t$
and $u$ is
\[
f^{\D_n}(t)=a_{rj}f^{\D_n}(u)^j+\dots
\]
Let $\lambda=v_X(f^{\D_n}(u))$. Then $v_X(f^{\D_n}(t))=\lambda j$,
and, applying the equalities from Hamburger-Noether expansions
of $\C$ and $\D$ from bottom to top, we deduce that
$v_X(f^\D(t_0))=\lambda v_X(f^\C(t_0))$ and
$v_X(f^\D(u_0))=\lambda v_X(f^\C(u_0))$. It follows that
$E_{\D}=\lambda E_{\C}$, whence $\lambda<l$ and $E_{\D}\ge E_{\C}$.

We may choose $t_0$ and $u_0$ so that the local equations
of the components of $R_{L/K,\CX}$ passing through $P_\C$
are within $\{t_0,u_0\}$. Then the same is true for
$R_{L/K,\CX_n}$, $O_n$ and $\{t,u\}$ respectively, where
$\CX_n$ is the $n$th monoidal transformation of $\CX$.
It follows $L=K(x)$, $x^l=t^qu^s\eps$, where $\eps$ is invertible in a
neighborhood of $O_n$. We obtain
\begin{align*}
n_\C=1 & \implies n_{\C_n}=1 \\
& \implies l\nmid q j + s = q v_X(f^{\C_n}(t)) + s v_X(f^{\C_n}(u))\\
& \implies l\nmid \lambda(q j + s)
 = q v_X(f^{\D_n}(t)) + s v_X(f^{\D_n}(u))\\
& \implies n_{\D_n}=1 \\
& \implies n_{\D}=1. \qed
\end{align*}
\end{pf}

\begin{lm}
\label{close}
Let $\C,\D$ be primitive arcs on $\CX$, at least one of them
being not regular. Assume that $(\C.\D)>(M_{\C}+1)E_{\C}E_{\D}$.
Then for any regular arc $\CF$ on $\CX$ we have
\[
(\D.\CF)=(\C.\CF)\cdot\frac{E_\D}{E_\C}.
\]
\end{lm}

\begin{pf}
Without loss of generality we may assume that
$P_\C=P_\D=P_{\CF}=:O$.
Let $t_0,u_0$ be a system of local parameters in $O$ such that
$t_0$ is a local equation of $R^{\CF}$. Then
$(\C.\CF)=v_X(f^\C(t_0))$, and $(\D.\CF)=v_X(f^\D(t_0))$.
As in the proof of Lemma \ref{l}, we obtain
$v_X(f^\D(t_0))=\lambda v_X(f^\C(t_0))$ and
$v_X(f^\D(u_0))=\lambda v_X(f^\C(u_0))$ for some $\lambda$.
It follows $E_\D=\lambda E_\C$.
\qed
\end{pf}

\begin{pr}
\label{lM}
Let $\CX=\Spec A$, $\CY=\Spec A_1$, where
$A$ and $A_1$ are as in Prop. \ref{tamelift}, $l\ne p$.
Let $\C$ be a primitive arc on $\CX$ with $[\C]\ne F_T$, and $\C'$ its lifting.
Then
\[
M_{\C'}\le \frac l2 (M_{\C}+1)E_{\C}^2.
\]
\end{pr}

\begin{pf}
If $n_\C=l$, we have $M_{\C'}\le M_{\C}$ by Lemma \ref{M}.
We may, therefore, assume that $n_\C=1$. Let $t,u$ be local
parameters of $A_1$ such that $v_X(f^{\C'}(t))>v_X(f^{\C'}(u))$,
namely, $t,u$ are just $T_1,U$ in a suitable order. Let
\begin{align*}
f^{\C'}(t) = & a_{01} f^{\C'}(u) + a_{02} f^{\C'}(u)^2 + \dots
 + a_{0h_0}f^{\C'}(u)^{h_0} + f^{\C'}(u)^{h_0}z_1, \\
f^{\C'}(u) = &\qquad a_{12} z_1^2 + \dots + a_{1h_1}z_1^{h_1} + z_1^{h_1}z_2, \\
& \dots \\
z_{r-2} = &\qquad a_{r-1,2} z_{r-1}^2  + \dots
        + a_{r-1,h_{r-1}}z_{r-1}^{h_{r-1}} + z_{r-1}^{h_{r-1}}z_r, \\
z_{r-1} = &\qquad a_{r2}z_r^2 + a_{r3}z_r^3+\dots
\end{align*}
be the Hamburger-Noether expansion of $\C'$ in the basis $t,u$. Define a new
arc $\D'$ by means of equalities:
\begin{align*}
f^{\D'}(t) = & a_{01} f^{\D'}(u) + a_{02} f^{\D'}(u)^2 + \dots
 + a_{0h_0}f^{\D'}(u)^{h_0} + f^{\D'}(u)^{h_0}y_1, \\
f^{\D'}(u) = &\qquad a_{12} y_1^2 + \dots + a_{1h_1}y_1^{h_1} + y_1^{h_1}y_2, \\
& \dots \\
y_{r-3} = &\qquad a_{r-2,2} y_{r-2}^2  + \dots
        + a_{r-2,h_{r-1}}y_{r-2}^{h_{r-2}} + y_{r-2}^{h_{r-2}}y_{r-1}, \\
y_{r-2} = &\qquad a_{r-1,2} y_{r-1}^2  + \dots
        + a_{r-1,h_{r-1}}y_{r-1}^{h_{r-1}} + y_{r-1}^{h_{r-1}+1},
        \\
y_{r-1} = & X.
\end{align*}
We see immediately that $\C'$ and $\D'$ have exactly
$h_0+\dots+h_{r-1}=M_{\C'}$ common infinitely near points,
whence $(\C'.\D')\ge 2M_{\C'}$. Next, we see that
$v_X(y_{r-1})<v_X(z_{r-1})$, $v_X(y_{r-1})<v_X(z_{r-1})$,
whence
\begin{equation}
\label{CuD}
\begin{aligned}
v_X(f^{\D'}(u)) & < v_X(f^{\C'}(u)), \\
v_X(f^{\D'}(t)) & < v_X(f^{\C'}(t)).
\end{aligned}
\end{equation}
It follows $E_{\D'}<E_{\C'}$.

By Remark \ref{L2}, $\D'$ is a lifting of a primitive arc $\D$
on $\CX$. By Prop. \ref{tamelift} it follows from (\ref{CuD})
that $E_{\D}/n_D<E_{\C}/n_\C=E_{\C}$.

Let first $n_\D=1$. Then $E_{\C}>E_{\D}$, in particular, $\C$ is not
regular.
Then by Prop. \ref{tamelift} we have
$(\C'.\D')=l(\C.\D)$. Assume that Proposition does not hold for
$\C$ and $\C'$, then
\[
(\C.\D)> (M_{\C}+1)E_{\C}^2>(M_{\C}+1)E_{\C}E_{\D}.
\]
It follows from Lemma \ref{l} that $E_{\D}\ge E_{\C}$, a
contradiction.

Let finally $n_\D=l$. Then by Prop. \ref{tamelift} we have
$(\C'.\D')=(\C.\D)$. We have $E_{\D}<l E_{\C}$, and by Lemma \ref{l} we
obtain
\[
(\C.\D)\le (M_{\C}+1)E_{\C}E_{\D} < l(M_{\C}+1)E_{\C}^2,
\]
and
\[
M_{\C'}< \frac l2 (M_{\C}+1)E_{\C}^2.  \qed
\]
\end{pf}

\section{Artin-Schreier extension}
\label{AS}

In this section we give an explicit shape of a uniform sufficient jet
order in the case of a cyclic extension of degree $p$.

\begin{pr}
\label{ASmain}
Let $\CX$ be a regular $k$-surface with the function field $K$, $L=K(x)$, where
$x^p-x=a\in K$.
Assume that $R=R_{L/K,\CX}$ is a simply normal crossing divisor.
Let $F_1,\dots,F_s$ be the components of $R$. Put
\[
m_i=-v_i(a),
\]
where $v_i$ is the valuation on $K$ associated with $F_i$,
$i=1,\dots,s$. Put
\[
R_0=R_0(\CX,L/K)=(m_1+1)F_1+\dots+(m_s+1)F_s.
\]
Let $\C,\D\in\bU_\CX$ be such that $P_\C=P_\D$, and
\begin{equation*}
(\C.\D)\ge (\C.R_0)E_{\D}+E_{\C}E_{\D}\max(M_{\C},M_{\D}).
\end{equation*}
Then $L/K$ has equal wild jumps at $\C$ and $\D$.
\end{pr}

\begin{pf}
Throughout the proof, we denote $O=P_\C=P_\D$.

Apply induction on $\max(M_{\C},M_{\D})$. Let $M_{\C}=M_{\D}=0$.
If none of the components of $R$ passes through $O$, we have
$w^{(i)}_\C=0=w^{(i)}_\D$ for any $i$.

If 1 component of $R$ passes through $O$, we may assume that this
is $F_1$. Let $t=0$ be a local equation of $F_1$ at $O$. Choose
any $u\in\CO_{\CX,O}$ such that $t,u$ are local parameters at $O$.

If there are 2 components
of $R$ through $O$, we may assume that these are $F_1$ and $F_2$.
Let $t=0$, $u=0$ be their local equations at $O$.

Since $F_i$ are components of $R$, each of the valuations $v_i$
must be totally or fiercely ramified in $L/K$; therefore, $m_i>0$.
We have $(\C.R_0)\ge m_1+1\ge2$, whence $N:=(\C.\D)\ge2$. Assume first that
$\C$ is not tangent to the curve locally defined as $u=0$.
Applying Weierstra{\ss} preparation theorem, we
can write $\CR^\C=k[[t,u]]/(f)$, $\CR^\D=k[[t,u]]/(g)$, where
\begin{align*}
f&=-t+\be_1 u + \be_2 u^2 +\dots, \\
g&=-t+\be'_1 u + \be'_2 u^2 +\dots, \\
\end{align*}
with $\be_i,\be_i'\in k$, $i=1,2,\dots$.
We have $\min\{i|\be_i\ne\be'_i\}=N$.

If $F_2$ passes through $O$, we have
\[
a=t^{-m_1}u^{-m_2}\eps(t,u),
\]
in the completion of $\CO_{\CX,O}$,
where $\eps\in k[[T,U]]$.
Now, $w^{(1)}_\C$ and  $w^{(1)}_\D$ are exactly the jump
of the Artin-Schreier equation
\begin{equation}
\label{A}
x^p-x=(\be_1 u+\be_2 u^2+\dots)^{-m_1}u^{-m_2}\eps(\be_1 u+\be_2
u^2+\dots,u),
\end{equation}
and that of
\begin{equation}
\label{B}
x^p-x=(\be'_1 u+\be'_2 u^2+\dots)^{-m_1}u^{-m_2}\eps(\be'_1 u+\be'_2
u^2+\dots,u).
\end{equation}
Let $j=(\C.F_1)$. Then $j=\min\{i|\be_i\ne0\}$. We have
\begin{align*}
&(\be_1 u+\be_2 u^2+\dots)^{-m_1}u^{-m_2}\equiv \\
&\qquad\equiv u^{-m_1j-m_2}(\be_j+\dots+\be_{j+m_1j+m_2-1}u^{m_1j+m_2-1})^{-m_1} \mod
k[[u]],
\end{align*}
and we see that the RHS of (\ref{A}) is congruent to that of
(\ref{B}) $\mod k[[u]]$, provided that $N\ge j+m_1j+m_2$.
This condition holds in view of $(\C.\D)\ge(\C.R_0)=(m_1+1)j+m_2+1$.
It follows $w^{(1)}_\C=w^{(1)}_\D$.

If $F_2$ does not pass through $O$, the argument is the same;
we have only to omit the factor $u^{-m_2}$ in all the formulae.

Finally, let $\C$ be tangent to the curve $u=0$. If $u=0$ is the
local equation of $F_2$, this case is reduced to the previous
one by changing the roles of $F_1$ and $F_2$. If there is no
component of $R$ passing through $O$ except for $F_1$, this
case can be also reduced to the previous one by a substitution
$u:=u+t$.

Next, let $\max(M_{\C},M_{\D})>0$. Let $\CX'\to\CX$ be the blowing up at $O$.
Let $\C',\D'$ be the strict transforms
of $\C$ and $\D$ respectively. Let $F_i'$ be the strict transform
of $F_i$, $i=1,\dots,s$, and let $E$ be the exceptional divisor.

Assume first that two components of $R$ (say, $F_1$ and $F_2$)
pass through $O$. The formula
\[
(\C.\D)\ge
((m_1+1)(\C.F_1)+(m_2+1)(\C.F_2))E_{\D}+E_{\C}E_{\D}\max(M_{\C},M_{\D})
\]
implies
\begin{align*}
(\C'.\D')+E_{\C}E_{\D}&\ge ((m_1+1)((\C'.F'_1)+E_{\C}) \\
&\qquad +(m_2+1)((\C'.F'_2)+E_{\C}))E_{\D} \\
&\quad+E_{\C}E_{\D}(\max(M_{\C'},M_{\D'})+1).
\end{align*}
and
\begin{align*}
(\C'.\D') &\ge ((m_1+1)(\C'.F'_1) +(m_2+1)(\C'.F'_2)+(m_1+m_2+1)(\C'.E))E_{\D} \\
&\quad +E_{\C}E_{\D}(\max(M_{\C'},M_{\D'})+1).
\end{align*}
Notice that $-v_E(a)\le m_1+m_2$, where $v_E$ is
the valuation associated with $E$. We obtain
\[
(\C'.\D')\ge (\C'.R_0(\CX',L/K))E_{\D'}+E_{\C'}E_{\D'}(\max(M_{\C'},M_{\D'})+1).
\]
Therefore, by the induction hypothesis, $L/K$ has equal
wild jumps at $\C'$ and $\D'$. Since $W_{\C,L/K}=W_{\C',L/K}$,
and $W_{\D,L/K}=W_{\D',L/K}$, $L/K$ has equal
wild jumps at $\C$ and $\D$ as well.

In the case when only one component of $R$ (say, $F_1$)
passes through $O$ the argument is similar.
\qed
\end{pf}

\begin{cor}
\label{AScor}
If $L/K$ is a cyclic extension of degree $p$,
Theorem \ref{main} holds for $\CX$ and $L/K$.
\end{cor}

\begin{pf}
Take
\[
\Delta_\CX(L/K,i)=i\max(m_1+1,\dots,m_s+1,i).
\qed
\]
\end{pf}

A slight modification of the above argument enables us
to prove the following analog of Lemma \ref{l}.
The assumptions about $\CX$ and $L/K$ are as in Prop. \ref{ASmain};
$n_\C$ denotes the number of non-equivalent liftings of $\C$.

\begin{lm}
\label{p}
Let $\C$ be a primitive arc on $\CX$ with $n_\C=1$.
Let $\D$ be a primitive arc on $\CX$ such that $E_{\D}<p E_{\C}$, and
\[
(\C.\D)>(\C.R_0(\CX,L/K))E_{\D}+E_{\C}E_{\D}M_{\C},
\]
where $R_0(\CX,L/K)$ is as in Prop. \ref{ASmain}.
Then

1. $n_\D=1$;

2. $E_{\D}\ge E_{\C}$.
\end{lm}

\begin{pf}
Let $n=M_{\C}$. Then $\C$ and $\D$ have at least $n$ common
infinitely near points $O_1,\dots,O_n$.
Choose any basis of local parameters at $O_0=P_\C=P_\D$
such that local equations of all components of $R_{L/K,\CX}$
passing through $O_0$ are within this basis. This choice determines the choice
of a distinguished basis of local parameters with the same property
at each of the $O_i$. Let $\CX_i$ be the blowing up
of $\CX_{i-1}$ at $O_{i-1}$, $\CX_0=\CX$.
Let $r$ be the number of lines in the Hamburger-Noether expansion
of $\C$. Considering Hamburger-Noether expansions
in the chosen basis, we see from $(\C.\D)>E_{\C}E_{\D}M_{\C}$
that the first $r-1$ lines in the expansions
of $\C$ and $\D$ are identical.

Let   $\C_n$ and $\D_n$ be the arcs which are the $n$th
strict transforms in $\CX_n$ of $\C$ and $\D$ respectively.
Since $n_\C=1$, we see that $w^{(1)}_{\C_n}=w^{(1)}_\C>0$. In
particular, at least one component of $R_{L/K,\CX_n}$ (say, $F_1$)
passes through $O_n$. Let $t,u$ be the distinguished basis at
$O_n$. Then we may assume without loss of generality that $t$
is the local equation of $F_1$, and $u$ is the local equation of
$F_2$, if there exists another component $F_2$ of $R_{L/K,\CX_n}$ passing through
$O_n$.

An argument similar to the step of induction in the previous
proof shows that
\begin{equation}
\label{bd}
(\C_n.\D_n)> (\C_n.R_0(\CX_n,L/K))E_{\D}.
\end{equation}

Consider the Hamburger-Noether expansion of the regular arc $\C_n$ in $t,u$:
\[
f^{\C_n}(t)=\be_1f^{\C_n}(u)+\be_2f^{\C_n}(u)^2+\dots
\]
Let
\[
j=\min\{i|\be_i\ne0\}=(\C_n.F_1).
\]
We see from (\ref{bd}) that $(\C_n.\D_n)>j E_{\C_n}E_{\D_n}$,
whence $\C_n$ and $\D_n$ have at least $j$ common infinitely
near points. In particular, the first line of the Hamburger-Noether
expansion of $\D_n$ in $t,u$ begins as
\[
f^{\D_n}(t)=\be_jf^{\D_n}(u)^j+\dots,
\]
and
\[
v_X(f^{\D_n}(t))=j v_X(f^{\D_n}(u)).
\]
 Looking at the first
$r-1$ lines of the Hamburger-Noether expansion of $\C$ and $\D$,
we conclude that $E_{\D}=E_{\C}v_X(f^{\D_n}(u))$, whence $E_{\D}\ge
E_{\C}$, and $v_X(f^{\D_n}(u))<p$.

Note that $R^{\C_n}=k[[T,U]]/(-T+B_0(U))$, where
\[
B_0(U)=\be_j U^j + \be_{j+1}U^{j+1}+\dots
\]
By the definition,
\[
(\C_n.\D_n)=v_X(-f^{\D_n}(t)+B_0(f^{\D_n}(u))),
\]
and we can write
\[
f^{\D_n}(t)=B_0(f^{\D_n}(u))+\delta,
\]
where $v_X(\delta)>(\C_n.R_0(\CX_n,L/K))E_{\D}$.

Assume that both $F_1$ and $F_2$ pass through $O_n$.
(The case when only $F_1$ passes through $O_n$ is similar.)
Then $L/K$ corresponds to an Artin-Schreier equation
\[
x^p-x=t^{-m_1}u^{-m_2}\eps(t,u),
\]
in the completion of $\CO_{\CX_n,O_n}$,
where $\eps\in k[[T,U]]$.
Then $w^{(1)}_{\C_n}$ and $w^{(1)}_{\D_n}$ are the ramification
jumps of the equations
\[
x^p-x=B_0^{-m_1}X^{-m_2}\eps(B_0,X),
\]
and
\[
x^p-x=f^{\D_n}(t)^{-m_1}f^{\D_n}(u)^{-m_2}\eps(f^{\D_n}(t),f^{\D_n}(u)),
\]
respectively. Note that
\begin{align*}
&
f^{\D_n}(t)^{-m_1}f^{\D_n}(u)^{-m_2}\eps(f^{\D_n}(t),f^{\D_n}(u))=
\\
& \qquad = (B_0(f^{\D_n}(u))+\delta)^{-m_1}f^{\D_n}(u)^{-m_2}\eps(B_0(f^{\D_n}(u))+\delta,f^{\D_n}(u))
\\
& \equiv
B_0(f^{\D_n}(u))^{-m_1}f^{\D_n}(u)^{-m_2}\eps(B_0(f^{\D_n}(u)),f^{\D_n}(u))\mod
k[[X]].
\end{align*}
It follows that $w^{(1)}_{\D_n}$ is equal to the ramification
jump of
\[
x^p-x= \Lambda(f^{\D_n(u)}),
\]
where $\Lambda=B_0^{-m_1}X^{-m_2}\eps(B_0,X)$.
Taking into account that $v_X(f^{\D_n(u)})<p$, we conclude that $w^{(1)}_{\D_n}$
is exactly $v_X(f^{\D_n(u)})$ times the ramification jump of
$x^p-x=\Lambda$, which is $w^{(1)}_{\C_n}$. In particular,
$w^{(1)}_{\D_n}>0$.
\qed
\end{pf}

Now we prove an analog of Prop. \ref{lM}.

\begin{pr}
\label{pM}
Let $\CX=\Spec A$, $\CX'=\Spec A_1$, where
$A$ and $A_1$ are as in Prop. \ref{tamelift}, $l=p$.
Let $\C$ be a primitive arc on $\CX$ with $[\C]\ne F_T$, and $\C'$ its lifting.
Then
\[
M_{\C'}<\frac p2(\C.R_0(\CX,L/K))E_{\C}+\frac p2 M_{\C}E_{\C}^2.
\]
\end{pr}

\begin{pf}
If $n_\C=p$, we have $M_{\C'}\le M_{\C}$ by Prop. \ref{tamelift}.
We may, therefore, assume that $n_\C=1$. Introduce $\D'$ and $\D$
exactly as in the proof of \ref{lM}. Again, we have
$E_{\D}/n_\D<E_{\C}$.

Let first $n_\D=1$. Then $E_{\C}>E_{\D}$.
By Lemma \ref{p}, we have
\[
(\C.\D) \le (\C.R_0(\CX,L/K))E_{\D}+M_{\C}E_{\C}E_{\D}.
\]
Next, by Prop. \ref{tamelift} we have $(\C'.\D')=p(\C.\D)$, whence
\begin{equation}
\label{fin}
M_{\C'}\le\frac12(\C'.\D')<\frac p2(\C.R_0(\CX,L/K))E_{\C}+\frac p2
M_{\C}E_{\C}^2.
\end{equation}

Let finally $n_\D=p$. Then by Prop. \ref{tamelift} we have
$(\C'.\D')=(\C.\D)$. On the other hand, $E_{\D}<p E_{\C}$, and Lemma
\ref{p} implies
\begin{align*}
(\C.\D) &\le (\C.R_0(\CX,L/K))E_{\D}+M_{\C}E_{\C}E_{\D} \\
&<p(\C.R_0(\CX,L/K))E_{\C}+p M_{\C}E_{\C}^2,
\end{align*}
and we have (\ref{fin}) again. \qed
\end{pf}

\section{Some reductions}

In this section we show that Theorem \ref{main} holds for $\CX$ if it
holds for some natural modifications of $\CX$.

Throughout this section $\CX$ is a surface over $k$, and
$L/K$ is a finite solvable Galois extension of the field
of functions on $\CX$.

\subsection*{Zariski covering}

\begin{pr}
\label{cov}
Let $\CX=\CX_1\cup\CX_2$, where $\CX_1$, $\CX_2$ are open subschemes
in $\CX$.
Then Theorem \ref{main} holds for $\CX$
and $L/K$ if it holds for $\CX_1$ and $L/K$ as well as for $\CX_2$ and $L/K$.
\end{pr}

\begin{rk} If, say, $\CX_1$ is not 2-dimensional, there are no arcs
on it. We agree that Theorem \ref{main}
is satisfied for $\CX_1$ in this case.
\end{rk}

\begin{pf}
It is sufficient to put
\[
\Delta_\CX(L/K,i)=\max(\Delta_{\CX_1}(L/K,i),\Delta_{\CX_2}(L/K,i))
\]
for any $i$.
If, say, $\CX_1$ is not 2-dimensional, we assume
$\Delta_{\CX_1}(L/K,i)=0$.
\qed
\end{pf}

\subsection*{Monoidal transformation}

Here we shall deduce the existence of a uniform sufficient jet order
for a surface $\CX$ from that for its monoidal transformation,
$\CX_1$.

\begin{pr}
\label{QT}
 Let $\CX_1\to\CX$ be a monoidal transformation.
Then Theorem \ref{main} holds for $\CX$
and $L/K$ if it holds for $\CX_1$ and $L/K$.
\end{pr}

\begin{pf}
Denote by $O\in\CX$ the closed point being blown up and by $E$
the exceptional curve. Let $F_1,\dots,F_n$ be all irreducible
components of $R_{L/K,\CX}$. Then, obviously, the irreducible
components of $R_{L/K,\CX_1}$ are $F^1_1,\dots,F^1_n$ and, possibly,
$E$, where $F_j^1$ is the strict transform of $F_j$, $j=1,\dots,n$.

Assume that Theorem \ref{main} holds
for $\CX_1$, and  $\Delta_{\CX_1}(L/K,-)$ is as in Theorem \ref{main}.
Assume first that $O$ lies on at least one of $F_i$, say, on $F_1$.
Introduce
\[
\Delta_\CX(L/K,i)=2\Delta_{\CX_1}(L/K,i)+i^2.
\]
Now let $\C,\D\in\bU_{\CX}$ be such that
\[
(\C.\D)\ge((\C.R)+\max(M_\C,M_\D))\Delta_\CX(L/K,\max(E_{\C},E_{\D})),
\]
and $P_\C=P_\D=:P$.
Let $\C'$ and $\D'$ be the strict transforms of $\C$ and $\D$
respectively.
If $P\ne O$, we have
\begin{align*}
(\C'.\D')&=(\C.\D) \\
& >((\C.R)+\max(M_\C,M_\D))\Delta_\CX(L/K,\max(E_{\C},E_{\D})) \\
&
\ge((\C.R_{L/K,\CX_1})+\max(M_{\C'},M_{\D'}))\Delta_{\CX_1}(L/K,\max(E_{\C'},E_{\D'})).
\end{align*}
Let $P=O$. Then it follows from (\ref{iqt}) that
\[
(\C'.\D')\ge((\C.R)+\max(M_\C,M_\D))\Delta_\CX(L/K,\max(E_{\C},E_{\D}))-E_{\C}E_{\D}.
\]
We have (e.~g., also by (\ref{iqt})) that $(\C.F_1)\ge E_{\C}$,
and $E_{\C}=(\C'.E)$ by Lemma \ref{EE}. Next,
\begin{align*}
&(\C.R)\Delta_\CX(L/K,\max(E_{\C},E_{\D})) \\
&\qquad \ge 2(\C.R)\Delta_{\CX_1}(L/K,\max(E_{\C},E_{\D})) + E_\C E_\D,
\end{align*}
whence
\begin{align*}
(\C'.\D') & > 2((\C.R)+\max(M_\C,M_\D))\Delta_{\CX_1}(L/K,\max(E_{\C},E_{\D}))
\\
& \ge ((\C.F_1)+\dots+(\C.F_n)+E_{\C}+\max(M_\C,M_\D))
\\
&\quad \times\Delta_{\CX_1}(L/K,\max(E_{\C},E_{\D}))
\\
& \ge((\C'.R_{L/K,\CX_1})+\max(M_{\C'},M_{\D'}))
\\
&\quad \times\Delta_{\CX_1}(L/K,\max(E_{\C'},E_{\D'})).
\end{align*}
Since Theorem \ref{main} holds for $\CX_1$ and $L/K$,
we conclude
\[
w^{(i)}_{\C}=w^{(i)}_{\C'}=w^{(i)}_{\D'}=w^{(i)}_{\D}
\]
for all $i=1,2,\dots$

Consider the remaining case, when $O$ lies on none of $F_i$.
In this case introduce
\[
\Delta_\CX(L/K,i)=\Delta_{\CX_1}(L/K,i).
\]
If $P\ne O$, we proceed as in the former case. Finally,
if $P=O$, both $\C$ and $\D$ do not meet the branch
locus, whence $w^{(i)}_{\C}=0=w^{(i)}_{\D}$ for $i=1,2,\dots$
\qed
\end{pf}


\subsection*{Unramified extension}

\begin{pr}
\label{URE}
 Let $\CX$ be a surface over $k$, $L/K$ a finite Galois extension
of the fraction field of $\CX$, $K'/K$ an unramified extension,
and $\CX'$ the normalization of $\CX$ in $K'$.
Then Theorem \ref{main} holds for $\CX$ and $L/K$
iff it holds for $\CX'$ and $K'L/K'$
(with $\Delta_\CX(L/K,-)=\Delta_{\CX'}(K'L/K',-)$).
\end{pr}

\begin{pf}
Denote by $f:\CX'\to\CX$ the normalization morphism.
Then $f^*R_{L/K,\CX}=R_{K'L/K',\CX'}$.
Let $\C$, $\D$ be arcs on $\CX'$ with $P_\C=P_\D=:P'$,
and let $P=f(P')$ be a regular point.
The ring $\CO_{\CX',P'}$ is unramified over $\CO_{\CX,P}$,
whence $\widehat{\CO_{\CX',P'}}=\widehat{\CO_{\CX,P}}$, and
\begin{align*}
(\C.\D) &=(f\circ\C.f\circ\D), \\
E_{\C} &=E_{f\circ\C}, \\
M_{\C} &=M_{f\circ\C}, \\
w^{(i)}_{\C}(K'L/K') &=w^{(i)}_{f\circ\C}(L/K), \\
w^{(i)}_{\D}(K'L/K') &=w^{(i)}_{f\circ\D}(L/K).
\end{align*}
Note that for any component $F$ of $R_{L/K,\CX}$, $f^*F$
is the sum of $[K':K]$ distinct prime divisors with no common
points.
Then $(f\circ\C.R_{L/K,\CX})=(\C.R_{K'L/K',\CX'})$
for any arc $\C$ on $\CX'$.

Let Theorem \ref{main} hold for $\CX$ and $L/K$.
Put
\[
\Delta_{\CX'}(K'L/K',i)=\Delta_{\CX}(L/K,i),
\]
and we see that
Theorem \ref{main} holds for $\CX'$ and $K'L/K'$.

Conversely, let Theorem \ref{main} be true for $\CX'$ and $L/K$.
Introduce
\[
\Delta_{\CX}(L/K,i)=\Delta_{\CX'}(K'L/K',i).
\]
Any arc $\C$ on $\CX$ can be written as $\C=f\circ\C'$,
where $\C'$ is an arc on $\CX'$.
We conclude that
Theorem \ref{main} holds for $\CX$ and $L/K$.
\qed
\end{pf}

\section{Proof of main theorem}

\begin{lm}
\label{tamemain}
Let $\CX$ be a surface over $k$ and $L/K$ a finite tame extension
of the fraction field of $\CX$. Then Theorem \ref{main} holds
for $\CX$ and $L/K$.
\end{lm}

\begin{pf}
In this case we have $w^{(i)}_\C=0$ for any $i$ and any suitable $\C$.
\qed
\end{pf}

A morphism of $k$-surfaces $\CX_2\to\CX_1$ is said to be \term{tame}
if it is dominant and proper and the corresponding extension
of function fields is tame. (Such morphisms are in fact
compositions of normalizations in tame extensions and birational
morphisms.)

\begin{lm}
\label{BCt}
Let $\CX_2\to\CX_1$ be a tame morphism, and let
$L/k(\CX_1)$ be a finite extension. Denote by $\CX_1'$
(resp., $\CX_2'$) the normalization of $\CX_1$
(resp., $\CX_2$) in $L$ (resp., $L k(\CX_2)$).
Then there exists a tame morphism $\CX'_2\to\CX'_1$
such that the corresponding field extension is $L k(\CX_2)/L$.
\end{lm}

\begin{pf}
Since $\CX_1'\to\CX_1$ is finite, any valuation
associated with a prime divisor on $\CX_1'$ s an extension
of a valuation
associated with a prime divisor on $\CX_1$. Therefore,
$L k(\CX_2)/L$ is tame with respect to $\CX_1'$.
It remains to note that the composition
$\CX_2'\to\CX_2\to\CX_1$ factors through $\CX_1'$.
\qed
\end{pf}

\begin{pr} Let $\CX$ be a surface over $k$ and $L/K$ a finite
solvable Galois extension of the fraction field of $\CX$.
Let $\CX'$ be the normalization of $\CX$ in a finite
solvable Galois extension $K'$ of $K$.
Assume that for any tame morphism of $k$-surfaces $\CX''\to\CX'$,
Theorem \ref{main} holds for $\CX''$ and $k(\CX'')L/k(\CX'')$.
Then Theorem \ref{main} holds for $\CX$ and $L/K$.
\end{pr}

\begin{pf} The proof consists of 3 steps.

\paragraph{1}
We apply induction on $[K':K]$.

Proposition is trivial for $[K':K]=1$.
Take $N>1$, and assume that Proposition holds for all extensions
of degree smaller than $N$.

Denote by $K_1/K$ any cyclic subextension in $K'/K$ of prime degree and
by $\CX_1$ the normalization of $\CX$ in $K_1$. By the induction hypothesis,
Theorem \ref{main} holds for $\CX_1$ and $K_1L/K_1$. In view of
Lemma \ref{BCt}, the same is true if one replaces $\CX_1$
with any $\CX_1''$ such that there exists
a tame morphism $\CX_1''\to\CX_1$.

Thus, we have reduced Proposition to the following case:
$K'/K$ is a cyclic extension of prime degree $l$.
We may also assume that $\CX$ is regular just removing
all singular points from it.

\paragraph{2}
Assume first that $l\ne p$.

First of all, we may assume without loss of generality that
\[
R_{L/K,\CX}\cup R_{K'/K,\CX}=R_{K'L/K,\CX}=:R
\]
is a simply normal
crossing divisor. Indeed, let $\CX_n$ be as in Proposition
\ref{NC}, and $\CX_n'$ be the normalization of $\CX_n$ in $K'$.
Then the proper morphism $\CX_n'\to\CX$ factors through a proper
birational morphism $\CX_n'\to\CX'$. If $\CX''\to\CX_n'$ is a
proper dominant morphism of $k$-surfaces such that $k(\CX'')=K''$
is a finite tame extension of $K'$ with respect to $\CX_n'$, then
$\CX''\to\CX'$ is also a proper dominant morphism of $k$-surfaces,
and $K''/K'$ is obviously tame with respect to $\CX'$. If
Proposition holds for $\CX_n$, it holds also  for $\CX$ by
Proposition \ref{QT}.

In view of \ref{cov},
we may assume that all the components of $R_{K'L/K,\CX}$ pass through
some $P\in\CX$, and $\CX=\Spec A$ is a sufficiently small
affine neighborhood of $P$. Let $t,u\in A$ be a system of local
parameters at $P$ which includes the local equations of all components of
$R_{K'L/K,\CX}$.
Taking $\CX$ small, we may also assume that $t$
and $u$ are prime elements in $A$, i.~e., the
Cartier divisor of $t$ (resp., $u$) is a prime divisor $F_t$
(resp., $F_u$).

If $R_{K'/K,\CX}$ is empty, Proposition follows immediately from
Prop. \ref{URE}. Next, consider the case when $R_{K'/K,\CX}$
consists of one component, and let $t\in A\subset\CO_{\CX,P}$ be the local
equation of this component. By Kummer theory, $K'=K(\root l\of a)$
for some $a\in K^*$. Let
\[
a=\eps t^{n_0}p_1^{n_1}\dots p_s^{n_s}
\]
be the canonical factorization of $a$ in $\CO_{\CX,P}$, where $\eps$
is invertible. The extension $K'/K$ is ramified at $(t)$ and unramified
at all $(p_i)$, whence $p$ divides (resp., does not divide) $n_i$ for
$i>0$ (resp., $i=0$). Choosing another generator of the subgroup
$\langle a\rangle(K^*)^p/(K^*)^p$ and another generator of the prime ideal
$(t)$, we may assume that $a=t$.

Let $t_1^l=t$. It is easy to see that all the singular points
of $\Spec A[t_1]$ lie over singular points of the curve $\Spec A/(t)$.
Note that this curve coincides locally with $R_{K'/K,\CX}$ and is therefore
regular at $P$.
If we replace $\CX$ with a smaller affine neighborhood of $P$, we may
assume that $\Spec A[t_1]$ is regular, whence $\CX'=\Spec A[t_1]$
is the normalization of $\CX$ in $K'$. We may also assume that $t_1$ is a prime element
of $A[t_1]$, i.~e., the
Cartier divisor of $t_1$ is a prime divisor $F_{t_1}$.

Denote by $f:\CX'\to\CX$ the normalization morphism.
Then
\[
f^*R\supset R_{K'L/K',\CX'}.
\]
 Let $Q$ be a point on $R$,
and let $Q'\in\CX'$ be a point above $Q$.
One of the components of $R$ passing through $Q$ has $t$ as
a local equation. If there is another component of $R$ through
$Q$, let $u_Q=u$ be its local equation (we have $Q=P$ in this case);
otherwise choose any $u_Q\in A$
such that $t,u_Q$ are local parameters at $Q$. It follows that
$t_1,u_Q$ are local parameters at $Q'$. Since $f^*R$ is locally
the zero locus of either $t_1$ or $t_1u_Q$, we conclude that
$R_{K'L/L,\CX'}$ is a simply normal crossing divisor.

Since Theorem \ref{main} holds for $\CX'$ and $K'L/K'$,
there exist some integers $\Delta_{\CX'}(K'L/K',i)$, $i=1,2,\dots$,
such that if
$\C',\D'\in\bU_{\CX'}$ are primitive arcs with $P_{\C'}=P_{\D'}$
and
\begin{equation}
\begin{aligned}
(\C'.\D')\ge & ((\C'.F_{t_1})+(\C'.F_u)+\max(M_{\C'},M_{\D'})) \\
&\quad \times \Delta_{\CX'}(K'L/K',\max(E_{\C},E_{\D})),
\end{aligned}
\label{eh}
\end{equation}
then $w_{\C'}^{(i)}(K'L/K')=w_{\D'}^{(i)}(K'L/K')$ for any $i$.
(If, say, $R_{K'L/K',\CX'}$ is just $F_{t_1}$, we omit the term $(\C'.F_u)$.)

Now introduce
\[
\Delta_\CX(L/K,i)=l^2i^2\Delta_{\CX'}(K'L/K',l i)
\]
for all positive integers $i$.

Let $\C,\D\in\bU_{\CX}$ be primitive arcs such that
\[
(\C.\D)\ge((\C.R)+\max(M_{\C},M_{\D}))\Delta_\CX(L/K,\max(E_{\C},E_{\D})).
\]
Let $Q=P_\C=P_\D$, $Q'\in\CX'$ be such that $Q=f(Q')$.
Let $Q\in F_t$. (We omit the trivial case $Q\notin F_t$.)
Localization and completion at $Q$ and $Q'$ respectively make $f$ into
the morphism $\Spec k[[T_1,U]]\to\Spec k[[T,U]]$ from Prop. \ref{tamelift}.
It follows from Prop. \ref{tamelift} that there exist primitive
arcs $\C',\D'$ which are liftings of $\C,\D$ respectively,
and $P_{\C'}=P_{\D'}=Q'$.
Fix $\D'$ and choose such lifting $\C'$ of $\C$ that $(\C'.\D')$
is maximal.
We have $(\C'.F_u)=l(\C.F_u)$. (This is Prop. \ref{tamelift}
if $Q=P$, and $0=0$ if $Q\ne P$.)

Applying further Prop. \ref{tamelift}, we obtain
\begin{align*}
n_\C n_\D(\C'.\D') & \ge l(\C.\D) \\
& \ge l ((\C.R)+\max(M_{\C},M_{\D}))\Delta_\CX(L/K,\max(E_{\C},E_{\D}))\\
& \ge (n_\C(\C'.R(K'L/K',\CX'))+l\max(M_{\C},M_{\D}))
  l^2\max(E_{\C},E_{\D})^2 \\
&\quad\times  \Delta_{\CX'}(K'L/K',l\max(E_{\C},E_{\D})) \\
& \ge (n_\C(\C'.R(K'L/K',\CX'))+l\max(M_{\C},M_{\D}))
  l^2\max(E_{\C},E_{\D})^2 \\
&\quad\times  \Delta_{\CX'}(K'L/K',\max(E_{\C'},E_{\D'})),
\end{align*}
whence
\begin{align*}
(\C'.\D')&\ge((\C'.R(K'L/K',\CX'))+\max(M_{\C},M_{\D})) \\
&\quad\times  l\max(E_{\C},E_{\D})^2\Delta_{\CX'}(K'L/K',\max(E_{\C'},E_{\D'}))\\
&\ge(l(\C'.R(K'L/K',\CX'))+l\max(M_{\C},M_{\D})\max(E_{\C},E_{\D})^2) \\
&\quad\times  \Delta_{\CX'}(K'L/K',\max(E_{\C'},E_{\D'})).
\end{align*}

By Prop. \ref{lM},
\[
M_{\C'}\le \frac l2(M_{\C}+1)E_{\C}^2 \le l M_{\C}E_{\C}^2,
\]
if $E_{\C}\ne1$, and $M_{\C'}\le\frac l2$ otherwise.
Therefore,
\[
\max(M_{\C'},M_{\D'}) \le l\max(M_{\C},M_{\D})\max(E_{\C},E_{\D})^2,
\]
if $E_{\C}\ne1$ or $E_{\D}\ne1$, and
\[
\max(M_{\C'},M_{\D'}) \le \frac l2,
\]
if $E_{\C}=E_{\D}=1$.
In both cases we conclude that
(\ref{eh}) is valid. It follows
\[
w_\C^{(i)}(L/K)=w_{\C'}^{(i)}(K'L/K')=w_{D'}^{(i)}(K'L/K')=w_\D^{(i)}(L/K)
\]
for any $i$.

It remains to consider the case when $R_{K'/K,\CX}$ consists of
two components.
Consider $K_1=K(t_1)$ and $K_2=K(t_1,u_1)$, where $t_1^l=t$,
$u_1^l=u$. Let $\CX_1$ and $\CX_2$ be the normalizations
of $\CX$ in $K_1$ and $K_2$ respectively. Then $K_2$ is a tame
extension of $K'$, whence Theorem \ref{main} is true
for $\CX_2$ and $K_2L/K_2$. Applying the already considered
case twice, we conclude that Theorem \ref{main} holds for $\CX_1$ and
$K_1L/K_1$, as well as for $\CX$ and $L/K$.

\paragraph{3}
Let now $l=p$.

Let $K'/K$ be given by Artin--Schreier equation
\[
x^p-x=a,
\]
$a\in K$.
An argument as in the case $l\ne p$ shows that we may replace $\CX$ with its
monoidal transformation.
After a suitable sequence of monoidal transformations,
we may assume that $R\cup\di a$ is a normal crossing divisor,
where $\di a$ is the divisor of $a$.
It follows that in the neighborhood
of any point $P$ we have $a=\eps t^iu^j$ in
$\CO_{\CX,P}$, where $\eps\in\CO_{\CX,P}^*$, and $t,u$ are local
parameters at $P$. Moreover, $L/K$ is unramified with respect
to any prime $\gp$ of $\CO_{\CX,P}$ unless $\gp=(t)$ or $\gp=(u)$.
Choosing a particular Artin-Schreier equation, we may require that either $i$
or $j$ is prime to $p$. Changing the respective local parameter,
we may assume $\eps=1$.

In view of \ref{cov},
we may assume that all the components of $R\cup\di a$ pass through
some $P\in\CX$, and $\CX=\Spec A$ is a sufficiently small
affine neighborhood of $P$. Taking $\CX$ small, we may also assume that the
divisor of $t$ (resp., $u$) is a prime divisor $F_t$
(resp., $F_u$).

Thus, we have
\[
\di a=iF_t+jF_u.
\]
We may also assume that at least one of $i,j$ is
negative. (Otherwise, $K'/K$ is unramified, and we argue as in the
case $l\ne p$.) We shall assume $i<0$. If $j>0$, we blow up at $P$
and obtain a scheme with
\[
\di a= iF_1+(i+j)E+jF_2,
\]
where $F_1$ is the strict transform of $F_t$, $E$ is the
exceptional divisor, $F_2$ is the strict transform of $F_u$.
Note that $|i+j|<\max(|i|,|j|)$. If $i+j>0$ (resp., $i+j<0$),
blow up at $F_1\cap E$ (resp., at $E\cap F_2$) and so on.
We obtain a scheme with
\[
\di a= n_0E_0+\dots+n_sE_s,
\]
where $E_\al$ does not meet $E_\be$ unless $\al=\be\pm1$,
$n_0,\dots,n_{r-1}<0$, $n_{r}=0$, $n_{r+1},\dots,n_s>0$
for some $r$. Locally we have the situation as before,
with $j\le0$.

We may assume that $(p,i)=1$. (Change the roles of $t$ and
$u$ if necessary. Note that we have already excluded the
case $p|i$, $p|j$.)
Next, making a tame base change of type $t=t_1^\al$, $u=u_1^\be$,
where $\al<p$, $\be<p$,
we may assume that $i\equiv1\mod p$, and either $j\equiv0\mod p$
or $j\equiv-1\mod p$. In the latter case,
blow up the point $P$. We get a scheme with
\[
\di a= iF_1+(i+j)E+jF_2,
\]
where $F_1$ (resp., $F_2$) is the strict transform of $F_t$ (resp., $F_u$),
$E$ is the exceptional divisor. We have $i+j\equiv0\mod p$.

Thus, it remains to consider the case
\[
x^p-x=t^{-mp+1}u^{-np},
\]
where $m>0$ and $n\ge0$. Put $t_1=t^mu^nx$. We have
\[
t_1^p-t^{m(p-1)}u^{n(p-1)}t_1=t,
\]
whence $t_1$ is integral over $A$, and $A[t_1]$ is regular.
Therefore, $\CX_1=\Spec A[t_1]$. We may also assume that the
divisor of $t_1$ (resp., of $u$) on $\CX_1$ is a prime divisor $F_{t_1}$
(resp., $F_u$).

Denote by $f:\CX'\to\CX$ the normalization morphism.
Then
\[
f^*R\supset R_{K'L/K',\CX'}.
\]
 Let $Q$ be a point on $R$,
and let $Q'\in\CX'$ be a point above $Q$.
Assume that $F_t$ passes through $Q$.
If there is another component of $R$ through
$Q$, let $u_Q=0$ be its local equation (we have $Q=P$ in this case);
otherwise choose any $u_Q\in A$
such that $t,u_Q$ are local parameters at $Q$. It follows that
$t_1,u_Q$ are local parameters at $Q'$.

Since Theorem \ref{main} holds for $\CX'$ and $K'L/K'$,
there exist some integers $\Delta_{\CX'}(K'L/K',i)$, $i=1,2,\dots$
such that if
$\C',D'\in\bU_{\CX'}$ are primitive arcs with $P_{\C'}=P_{\D'}$
and
\begin{equation}
\begin{aligned}
(\C'.\D')\ge & ((\C'.F_{t_1})+(\C'.F_u)+\max(M_{\C'},M_{\D'})) \\
&\quad\times \Delta_{\CX'}(K'L/K',\max{E_{\C},E_{\D}}),
\label{weh}
\end{aligned}
\end{equation}
then $w_{\C'}^{(i)}(K'L/K')=w_{\D'}^{(i)}(K'L/K')$ for any $i$.

Now introduce
\[
\Delta^\circ_\CX(L/K,i)=p^2 (i^2+\mu i^2)\Delta_{\CX'}(K'L/K',p i)
\]
for all positive integers $i$,
where $\mu=\max(m_1+1,\dots,m_s+1)$, $m_1,\dots,m_s$ are associated
with the extension $K'/K$ as in Proposition \ref{ASmain}.

Let $\C,\D\in\bU_{\CX}$ be primitive arcs such that
\begin{equation}
\label{cl}
(\C.\D)\ge((\C.R)+\max(M_{\C},M_{\D}))\Delta^\circ_\CX(L/K,\max(E_{\C},E_{\D})).
\end{equation}
Assume that the common closed point $Q$ of $\C$ and $\D$ lies on
$F_t$.

Similarly to the case $l\ne p$, we can construct liftings $\C'$ and $\D'$
of $\C$ and $\D$ such that
\begin{align*}
n_\C n_\D(C_1.D') & \ge p(\C.\D) \\
& \ge p ((\C.R)+\max(M_{\C},M_{\D}))\Delta_\CX(L/K,\max(E_{\C},E_{\D}))\\
& \ge p ((\C.R)+\max(M_{\C},M_{\D}))p^2(\max(E_{\C},E_{\D})^2+\mu\max(E_{\C},E_{\D})^2) \\
&\quad\times \Delta_{\CX'}(K'L/K',p\max(E_{\C},E_{\D})) \\
& \ge ((n_\C(\C'.R_{K'L/K',\CX'})+p\max(M_{\C},M_{\D}))
  p^2\max(E_{\C},E_{\D})^2 \\
&\qquad   +p(\C.R)p^2\mu\max(E_{\C},E_{\D})^2)  \\
&\quad\times  \Delta_{\CX'}(K'L/K',p\max(E_{\C},E_{\D})) \\
& \ge (n_\C(\C'.R_{K'L/K',\CX'})+p\max(M_{\C},M_{\D}))
  p^2\max(E_{\C},E_{\D})^2 \\
&\qquad   +p(\C.R)p^2\mu\max(E_{\C},E_{\D})^2)  \\
&\quad\times  \Delta_{\CX'}(K'L/K',\max(E_{\C'},E_{\D'})),
\end{align*}
whence
\begin{align*}
(\C'.\D')&\ge((\C'.R_{K'L/K',\CX'})+\max(M_{\C},M_{\D}))
  p(\max(E_{\C},E_{\D})^2 \\
&\qquad   +p(\C.R)\mu\max(E_{\C},E_{\D})^2)  \\
&\quad\times  \Delta_{\CX'}(K'L/K',\max(E_{\C'},E_{\D'})) \\
&\ge(p(\C'.R_{K'L/K',\CX'})+p\max(M_{\C},M_{\D})\max(E_{\C},E_{\D})^2 \\
&\qquad  +p(\C.R_0(\CX,K'/K))\max(E_{\C},E_{\D})^2)  \\
&\quad\times  \Delta_{\CX'}(K'L/K',\max(E_{\C'},E_{\D'})).
\end{align*}

By Prop. \ref{pM},
\[
M_{\C'}<\frac p2(\C.R_0(\CX,K'/K))E_{\C}+\frac p2 M_{\C}E_{\C}^2.
\]
It follows
\[
\max(M_{\C'},M_{\D'})<\frac p2((\C.R_0(\CX,K'/K))+ \max(M_{\C},M_{\D}))\max(E_{\C},E_{\D})^2.
\]
Indeed, this is trivial if both $\C$ and $\D$ are regular. In the opposite case,
it follows from (\ref{cl}) that
\[
(\C.\D)\ge (M_\C+1)\max(E_\C,E_\D)^2,
\]
and Lemma \ref{close} implies
\[
(\D.R_0(\CX,K'/K))=(\C.R_0(\CX,K'/K))\cdot\frac{E_\D}{E_\C}.
\]
From this, we obtain (\ref{eh}), whence
$w_{\C'}^{(i)}(K'L/K')=w_{\D'}^{(i)}(K'L/K')$ for any $i$.

In the case $Q\notin F_t$ we obtain that $w_{\C'}^{(i)}(K'L/K')=w_{\D'}^{(i)}(K'L/K')$
in a similar way.

Next, by Corollary \ref{AScor} we know that Theorem \ref{main} is true for $\CX$ and
$K'/K$. This means that there exists a non-decreasing sequence of positive numbers
$\Delta_\CX(K'/K,i)\ge i^2$, $i=1,2,\dots$
such that for any $\C,\D\in\bU_{\CX}$ with
\[
(\C.\D)\ge ((\C.R)+\max(M_{\C},M_{\D}))\Delta_\CX(K'/K,\max(E_{\C},E_{\D})),
\]
we have
\[
w^{(1)}_\C(K'/K)=w^{(1)}_\D(K'/K).
\]
Now let
\[
\Delta_\CX(L/K,i)=\max(\Delta^\circ_\CX(L/K,i),\Delta_\CX(K'/K,i)).
\]
Take any $\C,\D\in\bU_{\CX}$ with
\[
(\C.\D)\ge((\C.R)+\max(M_{\C},M_{\D}))\Delta_\CX(L/K,\max(E_{\C},E_{\D})).
\]
Then for suitable liftings $\C'$, $\D'$ of $\C$, $\D$
we have
\[
W_{\C,K'L/K}=W_{\C,K'/K}\circ W_{\C',K'L/K'}=
W_{\D,K'/K}\circ W_{\D',K'L/K'}=W_{\D,K'L/K}.
\]
Recall that $K'L/L$ is tame, whence
\[
W_{\C,L/K}=W_{\C,K'L/K}=W_{\D,K'L/K}=W_{\D,L/K},
\]
i.~e., $L/K$ has equal wild jumps at $\C$ and $\D$.
\qed
\end{pf}

Applying Proposition to the case $L=K'$ and using Lemma
\ref{tamemain}, we see that Theorem \ref{main} holds in full
generality.

\section{Towards equal ramification in the strong sense\label{strong}}

Let $\CX$ be as in the introduction,
$L/K$ a finite Galois extension of its function field,
and $\CY$ the normalization of $\CX$ in $L$.
Assume that $G=\Gal(L/K)$ is a $p$-group
with the following property: any cyclic subgroup of $G$
is normal in $G$.

\begin{rk}
There are some non-abelian groups with this property,
e.~g., the groups $\langle x,y|x^l=y^l,x^{l^2}=e,xy=y^{l+1}x\rangle$
for any prime $l$;
this was pointed out to me by V.~P.~Snaith.
\end{rk}

\begin{thm}
\label{T2}
There exists an effective divisor $R_0$ supported at $R=R_{L/K,\CX}$
such that for any $P\in R$ and any $C,C'\in U_P$
the condition $(C.C')_P\ge (C.R_0)_P$ implies that
$L/K$ is equally ramified with respect to $C$ and $C'$ at $P$ in the strong sense.
\end{thm}

\begin{pf}
In view of Corollary \ref{C1}, for any Galois subextension of
$M/K$ there exists an effective divisor $R_{0,M/K}$ supported
at $R_{M/K,\CX}\subset R_{L/K,\CX}$ such that $(C.C')_P\ge (C.R_{0,M/K})_P$
implies $W_{C,P,M/K}=W_{C',P,M/K}$.
Now fix any $R_0$ supported at $R=R_{L/K,\CX}$
such that $R_0\ge R_{M/K,\CX}$ for any $M/K$.

Let $C,C'\in U_P$, and $(C.C')_P\ge (C.R_0)_P$.
Let $C_L$ and $C'_L$ be any irreducible curves on $\CY$ over $C$ and
$C'$ respectively. Let $v$ be the valuation on $k(C)$ that corresponds
to $P$, and let $v_L$ be any extension of $v$ onto $k(C_L)$. Define $v'$ and $v'_L$
similarly. We have to prove that
$G_{j}(L_{C_L,v_L}/K_{C,v})=G_{j}(L_{C'_L,v'_L}/K_{C',v'})$
for any $j\ge-1$.

Choose any $g\in D_{C_L,v_L}$. It is sufficient to prove that
$g\in D_{C'_L,v'_L}$ and
\[
a_{C_L,v_L}(g)=a_{C'_L,v'_L}(g),
\]
where
\[
a_{C_L,v_L}(g)=\max\{j|g\in G_{j}(L_{C_L,v_L}/K_{C,v})\},
\]
and $a_{C'_L,v'_L}(g)$ is defined similarly.
Note that $D_{C_L,v_L}=G_{0}(L_{C_L,v_L}/K_{C,v})$ since the residue
field of the valuation on $K_{C,v}$ is the algebraically
closed field $k$. On the other hand, $L_{C_L,v_L}/K_{C,v}$
is a $p$-extension, whence
$G_{0}(L_{C_L,v_L}/K_{C,v})=G_{1}(L_{C_L,v_L}/K_{C,v})$.
Therefore, $a_{C_L,v_L}(g)$ (as well as $a_{C'_L,v'_L}(g)$)
is a positive integer.

Let $M$ be the fixed field of $\langle g\rangle$ inside $L/K$,
$\CZ$ the normalization of $\CX$ in $M$, $C_M$ (resp., $C'_M$)
the image of $C_L$ (resp., $C'_L$)  on $\CZ$,
$v_M$ (resp., $v'_M$)
the restriction of $v_L$ (resp., $v'_L$) onto the
function field of $C_M$ (resp., $C'_M$).
By the definition of $R_0$ we have $W_{C,P,L/K}=W_{C',P,L/K}$
and $W_{C,P,M/K}=W_{C',P,M/K}$, whence
\begin{equation}
\label{eqW}
\begin{aligned}
W_{L_{C_L,v_L}/M_{C_M,v_M}} &= W_{M_{C_M,v_M}/K_{C,v}}^{-1}\circ W_{L_{C_L,v_L}/K_{C,v}}
\\
&= W_{M_{C'_M,v'_M}/K_{C',v'}}^{-1}\circ W_{L_{C'_L,v'_L}/K_{C',v'}}
\\
&= W_{L_{C'_L,v'_L}/M_{C'_M,v'_M}}.
\end{aligned}
\end{equation}
Next,
\[
G_{j}(L_{C'_L,v'_L}/M_{C'_M,v'_M})=
G_{j}(L_{C'_L,v'_L}/K_{C',v'})\cap\Gal(L_{C'_L,v'_L}/M_{C'_M,v'_M}),
\]
whence
\begin{equation}
\label{a}
a_{C_L,v_L}(g)= \max\{j|g\in G_{j}(L_{C_L,v_L}/M_{C_M,v_M})\}.
\end{equation}
Let $p^r$ be the order of $g$ in $G$.
The condition $g\in D_{C_L,v_L}$ is equivalent to
$|G_1(L_{C_L,v_L}/M_{C_M,v_M})|=p^r$ and to
$w^{(r)}(L_{C_L,v_L}/M_{C_M,v_M})>0$.
Note that $|G_1(L_{C_L,v_L}/M_{C_M,v_M})|^{-1}$
is the minimal slope of the graph of
$W_{L_{C_L,v_L}/M_{C_M,v_M}}$. We obtain from (\ref{eqW})
that $g\in D_{C'_L,v'_L}$.

It is clear from (\ref{a}) that
\[a_{C_L,v_L}(g)=w^{(r)}(L_{C_L,v_L}/M_{C_M,v_M}).
\]
Since $w^{(r)}(L_{C_L,v_L}/M_{C_M,v_M})$ is just
the abscissa of the end of the first segment in the graph
of $W_{L_{C_L,v_L}/M_{C_M,v_M}}$, we conclude from (\ref{eqW})
that $a_{C_L,v_L}(g)=a_{C'_L,v'_L}(g)$.
\qed
\end{pf}

\begin{cor}
Let $G=\Gal(L/K)$ be a $p$-group
with the following property: any cyclic subgroup of $G$
is normal in $G$.
Then the answer to Question \ref{eusu-sw} is affirmative.
\end{cor}

\begin{pf}
See Remark \ref{suffstrong}.
\qed
\end{pf}


\bigskip

Department of Mathematics and Mechanics, St. Petersburg University

Bibliotechnaya pl., 2, Staryj Petergof

198904 St. Petersburg Russia

\medskip

E-mail: {\tt Igor.Zhukov@mail.ru, igor\_zh@hotmail.com}

\end{document}